\documentclass[12pt,leqno,color,amssymb]{amsart}

\usepackage{amssymb, amsmath,latexsym, amsfonts}

  \def\emty{\emptyset}
  \def\cbar{\widehat{\C}}
  
  \def\stab{\mbox{\rm Stab}\,}
  \def\diam{\mbox{\rm diam}}
  
  \def\dist{\mbox{\rm dist}}
  \def\rag{\rangle}
  \def\lag{\langle}
  
  \newcommand\ben{\begin{enumerate}}
  \newcommand\een{\end{enumerate}}
  \newcommand\bit{\begin{itemize}}
  \newcommand\eit{\end{itemize}}

  \def\CCC{{\mathcal C}}

  \def\FFF{{\mathcal F}}
  \def\GGG{{\mathcal G}}
  \def\HHH{{\mathcal H}}

  \def\PPP{{\mathcal P}}
  \def\QQQ{{\mathcal Q}}

  \def\UUU{{\mathcal U}}
  \def\VVV{{\mathcal V}}
  
  \def\YYY{{\mathcal Y}}

  \def\th{theorem }
  
  \def\Th{Theorem }

  \def\homeo{homeomorphism }

  \def\Proof{{\noindent\sc Proof.} }

  \def\al{\alpha}
  \def\be{\beta}
  \def\te{\theta}
  
  \def\g{\gamma}
  \def\G{\Gamma}
  
  \def\Si{\Sigma}
  \def\vp{\varphi}
  \def\ep{\varepsilon}
  \def\la{\lambda}
  \def\La{\Lambda}
  \def\De{\Delta}
  \def\de{\delta}
  
  \def\oo{\mbox{$\Omega$} }

  \def\R{\mbox{$\mathbb R$}}
  
  \def\C{\mbox{$\mathbb C$}}
  
  \def\Z{\mbox{$\mathbb Z$}}

  \def\SS{\mbox{$\mathbb S$}}
  \def\HH{\mbox{$\mathbb H$}}
  \def\DD{\mbox{$\mathbb D$}}
  \def\PP{\mbox{$\mathbb P$}}

  \def\heta{\widehat{\eta}}

  \def\endp{\hspace*{\fill}$\rule{.55em}{.55em}$ \smallskip}

      \newtheorem{newthm}{Theorem}
      \newtheorem{theorem}{Theorem}[section]
      \newtheorem{lemma}[theorem]{Lemma}
      \newtheorem{proposition}[theorem]{Proposition}
      \newtheorem{corollary}[theorem]{Corollary}
      \newtheorem{defthm}[theorem]{Definition et \th}
      \newtheorem{defn}[theorem]{Definition}
      \newcommand{\REFEQN}[1] { \begin{equation}\label{#1} }
      \newcommand{\ENDEQN}{\end{equation}}
      \newcommand{\REFTHM}[1] { \begin{theorem}\label{#1} }
      \newcommand{\ENDTHM}{\end{theorem}}
      \newcommand{\REFNTH}[1] { \begin{newthm}\label{#1} }
      \newcommand{\ENDNTH}{\end{newthm}}
      \newcommand{\REFPROP}[1]{\begin{proposition}\label{#1} }
      \newcommand{\ENDPROP}{\end{proposition} }
      \newcommand{\REFLEM}[1]{\begin{lemma}\label{#1} }
      \newcommand{\ENDLEM}{\end{lemma} }
      \newcommand{\REFCOR}[1]{\begin{corollary}\label{#1} }
      \newcommand{\ENDCOR}{\end{corollary} }
      \newcommand{\REFDEFTHM}[1] { \begin{defthm}\label{#1} }
      \newcommand{\ENDDEFTHM}{\end{defthm}}

\numberwithin{equation}{section}
\newcommand{\wLa}{\widehat{\La}}
\renewcommand{\hat}{\widehat}

\newcommand{\confdim}{\mbox{\rm confdim}}

\newtheorem{rmk}[theorem]{Remark}
\newtheorem{conj}[theorem]{Conjecture}

\title{Hyperbolic groups with planar boundaries
}
\author{Peter Ha\"{\i}ssinsky}
\address{Universit\'e Paul Sabatier,
Institut de Math\'ematiques de Toulouse,
118 route de Narbonne,
31062 Toulouse Cedex 9, France}
\email{phaissin@math.univ-toulouse.fr}

\date{\today}

\begin{document}



\begin{abstract} We prove that the class of convex-cocompact Kleinian groups is quasi-isometrically rigid.
We also establish that a word hyperbolic group with a planar boundary 
different from the sphere is virtually a convex-cocompact 
Kleinian group
provided that its boundary has Ahlfors regular conformal dimension
strictly less than $2$ or if it acts geometrically on a CAT(0) cube complex.\end{abstract}

\maketitle

\section{Introduction} 

A {\it Kleinian group}  is a discrete subgroup of $\PP SL_2(\C)$ which we view as acting
both on hyperbolic $3$-space $\HH^3$ via orientation-preserving isometries and on the Riemann
sphere $\cbar$ via M\"obius transformations. 
Since  Poincar\'e introduced them 
for solving differential
equations with algebraic coefficients at the end of the nineteenth century \cite{poincare:memoire:kleineen}, Kleinian groups have continuously drawn a lot of attention,
playing a prominent role in
complex analysis, conformal dynamical systems, hyperbolic geometry,  Teichm\"uller theory and low dimensional topology, see for instance \cite{klein:neue:func:theorie, bers:teichmueller:bams,
sullivan:bams, thurston:bams}. The subclass of
 {\it convex-cocompact Kleinian groups} is particularly relevant for the topology of $3$-dimensional manifolds, cf. Corollary \ref{cor:thurston}. 
These are finitely generated Kleinian groups $G$ for which there is a convex subset $\CCC\subset \HH^3$ invariant under $G$ such that $\CCC/G$ is compact. In particular, a {\it cocompact} Kleinian group is convex-cocompact with $\CCC=\HH^3$.

In the light of their prominent appearance in different areas of mathematics, characterizing this class of convex-cocompact Kleinian groups becomes a natural and interesting problem. The central purpose of this paper is to offer an array of such characterizations.

\subsection{Characterizations of convex-cocompact Kleinian groups}
Most of the characterizations proposed below will hold up to a finite index subgroup; in other words, we will show that groups are {\it virtually convex-compact Kleinian groups}.
The main 
reason for considering characterizations up to finite index comes from an example provided by Kapovich and Kleiner 
in \cite[\S\,8]{kapovich:kleiner:lowdim} which shows that having to pass to a finite index subgroup cannot be generally avoided without any further assumption.

\subsubsection{Characterization up to quasi-isometry}
From the point of view of geometric group theory, one tries to understand the properties of a group by studying the different
actions it admits on metric spaces. To start with, a finitely generated group acts by left-translation on the Cayley graph $X$ associated to any of its finite generating sets. If 
such a graph is equipped with the length metric which makes each edge isometric to the segment $[0,1]$, then  the action of $G$ becomes 
{\it geometric}: the group $G$ acts by isometries (the action is distance-preserving), properly discontinuously
(for any compact subsets $K$ and $L$ of $X$, at most finitely many elements $g$ of $G$ will satisfy $g(K)\cap L \ne\emty$) and 
cocompactly (the orbit space $X/G$ is compact). 

The classification of finitely generated groups up to {\it quasi-isometry} is a central issue:
a quasi-isometry between metric spaces $X$ and $Y$ is a map  $\vp:X\to Y$ such that there are constants $\la>1$ and $c>0$ such that:
\bit
\item (quasi-isometric embedding) for all $x,x'\in X$, the two inequalities $$\frac{1}{\la}d_X(x,x') - c \le d_Y(\vp(x),\vp(x'))\le \la d_X(x,y)+ c$$ hold and
\item the $c$-neighborhood of the image $f(X)$ covers $Y$.
\eit
This defines in fact an equivalence relation on metric spaces. Note that any two locally finite Cayley graphs of the same group are  quasi-isometric, which enables us to discuss 
the quasi-isometry class of a finitely generated group (through the class of its locally finite Cayley graphs). More generally,  
 \v{S}varc-Milnor's lemma asserts that there is  only one geometric action 
of a group on a geodesic metric space up to quasi-isometry \cite[Prop.\,3.19]{ghys:delaharpe:groupes}.

Our first result says that a group $G$ is virtually a convex-cocompact Kleinian group if it looks like one:

\begin{theorem}[quasi-isometric rigidity]\label{thm:mmain} The class of convex-cocompact Kleinian groups is quasi-isometrically rigid. More precisely,
 a finitely generated group quasi-isometric to a convex-cocompact Kleinian group  contains a finite
index subgroup isomorphic to a  (possibly different) convex-cocompact Kleinian group.
\end{theorem}

When $G$ is quasi-isometric to a cocompact Kleinian group, the theorem follows from \cite{cannon:cooper}. 
Theorem \ref{thm:mmain} covers two trivial cases which define the class of {\it elementary groups}: a group is elementary if it  is  finite or if it contains a finite index subgroup isomorphic to  $\Z$.  
Note that Ker\'ekj\'art\'o proved that compact groups of homeomorphisms of the sphere are always conjugate to subgroups of $O_3(\R)$ \cite{kerekjarto:acta,kolev:sousgroupe_compact}.
General background on quasi-isometric rigidity includes  \cite{drutu:qirigidity, kapovich:qirigidity}.
Actually a stronger form holds:

\begin{theorem}\label{thm:cmain} 
A finitely generated group admitting a quasi-isometric embedding into $\HH^3$
contains a finite index subgroup isomorphic to a  convex-cocompact Kleinian group.\end{theorem}

\begin{rmk}\label{rmk:gdcboundary} If a finitely generated group $G$ is virtually isomorphic to a convex-cocompact Kleinian group, 
then it is 
quasi-isometric to a convex subset of $\HH^3$ 
with geodesic boundary iff its  boundary is homeomorphic to
 the Sierpi\'nski carpet (see \S\,\ref{subsec:planaraction} for a proof).\end{rmk}


\subsubsection{Dynamics on the boundary}
The action of a Kleinian group $G$ on the Riemann sphere partitions $\cbar$ into the {\it domain of discontinuity} 
$\oo_G$, which is the largest set of $\cbar$ on which $G$ acts properly discontinuously, and the {\it limit
set} $\La_G$, which is the minimal $G$-invariant compact subset of $\cbar$. 
The action of a convex-cocompact Kleinian group on its limit set is that of a {\it uniform convergence group action},  i.e., the diagonal action on the set
of distinct triples is properly discontinuous (convergence action) and cocompact (uniform). It follows from Bowditch's topological characterization that $G$ is word hyperbolic in the sense
of Gromov and that $\La_G$ is a model for its boundary at infinity \cite{bowditch:characterization}. 
More generally, Bowditch's work enables us to define a {\it word hyperbolic group} $G$ as a group acting by homeomorphisms on a metrizable compact space $Z$ in
such a way that its diagonal action on the set of distinct triples is properly discontinuous and cocompact.
In the theory, the topological space $Z$ is a model for {\it the boundary at infinity $\partial G$ of $G$}. This characterization shows that this notion of hyperbolicity
is purely topological i.e., it only depends on the topological dynamics of the group on its boundary.
Background on hyperbolic groups include \cite{gromov:hyperbolic,ghys:delaharpe:groupes,
kapovitch:benakli}.


A general principle asserts that a word hyperbolic group is determined  by its boundary. More precisely, Paulin proved that the quasi-isometry class of a word hyperbolic group
is determined by its boundary equipped with
its quasiconformal structure \cite{paulin:determined}.  Let us recall some definitions.
A homeomorphism $h: X \to Y$ between metric spaces is called {\em quasisymmetric} provided 
there exists a homeomorphism $\eta: [0,\infty) \to [0,\infty)$ such that 
$d_X(x,a) \leq td_X(x,b)$ implies $d_Y(f(x),f(a)) \leq \eta(t)d_Y(f(x),f(b))$ for all triples of points 
$x, a, b \in X$ and all $t \geq 0$  \cite{tukia:vaisala:qs}. The boundary $\partial G$ of a hyperbolic group $G$ 
is endowed with a {\em  conformal
gauge} $\GGG(G)$ i.e., a family of metrics which are pairwise quasisymmetrically equivalent. These metrics are
exactly those metrics compatible with the topology of $\partial G$ for which the action of $G$ is {\it uniformly quasi-M\"obius} \cite{vam}: 
there exists a \homeo $\te:\R_+\to\R_+$
such that, for any distinct points $x_1,x_2,x_3,x_4\in \partial G$, for any $g\in G$,
$$[g(x_1):g(x_2):g(x_3):g(x_4)]\le \te([x_1:x_2:x_3:x_4])$$
where $$[x_1:x_2:x_3:x_4]= \frac{|x_1-x_2|\cdot |x_3-x_4|}{|x_1-x_3|\cdot |x_2- x_4|}\,.$$
Paulin's result says that two non-elementary word hyperbolic groups are quasi-isometric if and only if there is homeomorphism between their boundaries
identifying their conformal gauges.

%


The aforementioned quasi-isometric rigidity Theorem \ref{thm:mmain} is established using the boundary at infinity by means of
Theorem \ref{cor:main1}.

\REFTHM{cor:main1} Let $G$ be a non-elementary hyperbolic group and
let $\partial G$ be endowed with a metric from its conformal gauge. If there exists a quasisymmetric
embedding of $\partial G$ into $\cbar$, then $G$ is virtually isomorphic to a convex-cocompact Kleinian group.\ENDTHM

\begin{rmk}\label{rmk:qivsqs} {\rm Knowing that $G$ is hyperbolic, the assumption
that $G$ admits a quasi-isometric embedding into $\HH^3$ or that
its boundary embeds quasisymmetrically into $\cbar$ are equivalent,
see \cite[Thm.\,7.4, Thm.\,8.2 and Prop.\,10.1]{bonk:schramm:embeddings}}.\end{rmk}

In low dimension, it is expected that  topology should determine geometry. 
This is the case for the uniformization of compact surfaces and three-manifolds, which only depends on
their fundamental groups. As was mentioned above, a convex-cocompact Kleinian group is word hyperbolic and its limit set is {\it planar} i.e., 
it embeds continuously into  the two-sphere $S^2$.
Thus, in the present setting, this general principle that the geometry is determined by topological data has the following interpretation:

\begin{conj}\label{conj:main} A word hyperbolic group with planar boundary 
contains a finite index subgroup isomorphic to a convex-cocompact Kleinian group.
\end{conj}

This would imply that we may drop the knowledge of the  conformal gauge in Paulin's characterization of a word hyperbolic 
group when its boundary is planar.
The aim  of this paper
is to provide supporting evidence to this picture (Conjecture \ref{conj:main}).
It  is known to hold when the boundary is homeomorphic to a simple closed curve by the works
\cite{casson:jungreis, gabai:S1}, characterizing more generally Fuchsian groups as convergence groups acting on the unit circle. 
Besides, Theorem \ref{cor:main1}  provides us with a  positive answer for a weaker and somewhat intermediate conjecture, and it enables us to restate Conjecture \ref{conj:main}
 in analytic terms:

\begin{conj} If the boundary of a word hyperbolic group is planar, then it admits a quasisymmetric embedding in the
Riemann sphere sphere $\cbar$, when equipped with a metric of its conformal gauge.\end{conj}

The following two well-known conjectures can be derived from the above problem by specifying the boundary of the group:

\begin{conj} Let $G$ be a word hyperbolic group.
\bit
\item {\bf (Cannon, \cite[Conjecture 11.34]{cannon:conj})} If $\partial G$ is homeomorphic to $S^2$, then $G$ contains a finite index subgroup isomorphic to
a cocompact Kleinian group.
\item {\bf (Kapovich and Kleiner,  \cite[Conjecture 6]{kapovich:kleiner:lowdim})}  If $\partial G$ is homeomorphic to the Sierpi\'nski carpet, then $G$ contains a finite index subgroup isomorphic to
the fundamental group of a compact hyperbolic $3$-manifold with non-empty totally geodesic boundary.
\eit
\end{conj}


Recall that the Sierpi\'nski carpet is the metric space obtained by starting 
with the unit square, subdividing it into nine squares, removing the middle square, repeating this procedure
ad infinitum with the remaining squares, 
and taking the decreasing intersection.  We shall define a 
{\it carpet group} to be a hyperbolic group with a boundary homeomorphic to the 
Sierpi\'nski carpet.  Note that  Theorem \ref{cor:main1}  for carpet groups follows from Bonk's work  \cite{bonk:sierpinski}.
In that case, the conclusion
can be strengthened \cite{kapovich:kleiner:lowdim, frigerio:commensurability, bonk:kleiner:merenkov}, see \S\,\ref{sec:qirigidity} for more details.

\subsubsection{Ahlfors regular conformal dimension}
Our further evidence supporting Conjecture \ref{conj:main} is based on the 
 so-called {\em Ahlfors regular conformal dimension}
$\confdim_{AR}G$ of $G$.  A metric space $X$ is  {\em Ahlfors regular} 
if there is a  Radon measure $\mu$ such that
for any $x\in X$ and $r\in (0,\diam\,X]$,
$\mu(B(x,r))\asymp r^Q$ holds for some given $Q>0$ \cite{mattila:gmt}. 
The measure $\mu$ is equivalent to the Hausdorff measure of $X$ of dimension $Q$.
Given a word hyperbolic group $G$, the subset of Ahlfors regular metrics $\GGG_{AR}(G)$ in $\GGG(G)$
is non-empty  \cite{coornaert:patterson-sullivan} and defines the Ahlfors-regular conformal gauge of the group (or of
its boundary).
The  {\em Ahlfors regular conformal dimension}
$\confdim_{AR}G$ of $G$ is defined as 
the infimum over $\GGG_{AR}(G)$ of every
dimension $Q$  \cite{mackay:tyson:confdim, carrasco:confgauge, ph:bbki}. This is a numerical invariant of the quasi-isometry class of $G$ by \cite{paulin:determined},
and is also a topological invariant of the action of the  group. It is a refinement due to Bourdon and Pajot \cite{bourdon:pajot:besov}  of 
Pansu's notion of conformal dimension \cite{pansu:confdim}.


We will prove that this numerical invariant suffices to characterize convex-cocompact Kleinian groups.

\REFTHM{thm:main1} Let $G$ be a non-elementary hyperbolic group with planar boundary non-homeomorphic to
the sphere. Then  $G$ is virtually isomorphic to a convex-cocompact Kleinian group if and only if $\confdim_{AR} (G) <2$.
\ENDTHM
The necessity part of Theorem \ref{thm:main1} is due to Sullivan, see
\cite{sullivan:ihes50}. For carpet groups, the result was previously established by Bonk and Kleiner, see Corollary \ref{cor:main} and Remark \ref{rmk:bonkleiner}  for
more details. 
When the boundary is the whole sphere,  Bonk and Kleiner also proved that the group is Kleinian
if and only if there is a distance of minimal dimension in its Ahlfors-regular conformal gauge \cite{bonk:kleiner:conf_dim}.

\subsubsection{Groups acting on CAT(0) cube complexes}
In general, the conformal dimension is difficult to estimate, see nonetheless \cite{carrasco:gafa} for a criterion ensuring dimension one
for some word hyperbolic groups. One way to overcome this difficulty is to use an induction argument for groups admitting a quasiconvex hierarchy,
as for the fundamental groups of Haken manifolds, cf.  \S\,\ref{sec:lowtop}. Recent results of Wise and Agol \cite{wise:qcvxh, agol:virtualHaken} enable us to
express this property in terms of  groups acting on CAT(0) cube complexes (\S\,\ref{sec:cccat}):


\REFTHM{thm:main2} Let $G$ be a hyperbolic group with planar boundary. Then $G$ is virtually isomorphic to a convex-cocompact
Kleinian group if and only if  $G$ acts geometrically and cellularly on a CAT(0) cube complex.\ENDTHM

Here, we admit groups with boundary homeomorphic to the whole sphere, recovering Markovic's criterion for the  Cannon conjecture:

\begin{corollary}[Markovic]\label{cor:main2} Let $G$ be hyperbolic group with boundary homeomorphic to the sphere and
which has a faithful and orientation-preserving action on its boundary. Then $G$ is isomorphic to a cocompact Kleinian group
if and only if $G$ acts geometrically and cellularly on a CAT(0) cube complex.\end{corollary}

\begin{rmk} {\em The original proof of the sufficiency of Corollary \ref{cor:main2} consists in extending the action on the two-sphere to the unit ball as a free
convergence action so that the quotient is a Haken manifold, concluding then with Thurston's uniformization Theorem \ref{thm:thurston}  \cite{markovic:cannon}.  Here, we use the action on the cube complex 
to split the group in order to obtain groups  with one-dimensional boundary and so that  our previous results can be applied.}
\end{rmk}

\subsubsection{Reduction to carpet groups}
In some cases, we can obtain a hierarchy based on the accessibility of word hyperbolic groups over  elementary
subgroups \cite{delzant:potyagailo:access,louder:touikan}. 

\REFTHM{thm:main0} Let $G$ be a non-elementary hyperbolic group with planar boundary non-homeomorphic to
the sphere and with no elements
of order two. The following conditions are equivalent.
\ben
\item $G$ is virtually isomorphic to a convex-cocompact Kleinian group.
\item  $\confdim_{AR} (G) <2$.
\item $\confdim_{AR} (H) <2$ for all quasiconvex carpet subgroups $H$ of $G$.
\een
\ENDTHM
We refer to Section \ref{sec:hyp} for the notion of quasiconvexity.
The implication (2) $\Rightarrow$ (3) follows for instance from \cite[Prop.\,2.2.11]{mackay:tyson:confdim}.

The following statement is an immediate consequence of Theorem \ref{thm:main0}.

\REFCOR{cor:main0} Every hyperbolic group $G$ with a one-dimensional planar boundary and no elements of order two
is virtually Kleinian if and only if every carpet group is virtually Kleinian. In particular, 
if $G$ has no carpet subgroup, then $G$ is virtually Kleinian.\ENDCOR

Theorem \ref{thm:main0} and Corollary \ref{cor:main0} reduce the dynamical characterization 
of (convex) cocompact Kleinian groups (Conjecture \ref{conj:main}) to both
Cannon and Kapovich-Kleiner conjectures. Let us recall that 
the Kapovich-Kleiner conjecture is implied by the Cannon conjecture \cite[Thm.\,5, Cor.\,13]{kapovich:kleiner:lowdim}.
Thus, the dynamical characterization of 
convex-cocompact Kleinian groups with no $2$-torsion
would follow from the Cannon conjecture. On the other hand, 
the Cannon conjecture 
does not imply Theorems  \ref{thm:main1}, \ref{thm:main2} nor \ref{thm:main0} since it is not known that a word hyperbolic group with a planar boundary is virtually a quasiconvex
subgroup of  a word hyperbolic group with the sphere as boundary.

\subsubsection{Planar actions}
All these results will rely on a particular case for which we know that 
the action is already planar in the following
sense.
We shall say that a hyperbolic group 
has a {\it planar action} if its boundary admits a topological embedding into
the two-sphere in such a way that the action of every element of the group can be extended to a global homeomorphism. 
Recall that a word hyperbolic is one-ended if and only if its boundary is non-empty and  connected.

\REFTHM{thm:main} Let $G$ be a one-ended  hyperbolic group with a planar action and with
boundary  non-homeomorphic to the sphere.
Then $G$ is virtually isomorphic to a convex-cocompact Kleinian group if and only if $\confdim_{AR} G< 2$.
\ENDTHM

\begin{rmk} If, in the above theorem, the action is faithful and orientation-preserving on its boundary, 
then the proof shows that
the group is isomorphic to a convex-cocompact Kleinian group; see also Corollary \ref{cor:planaraction} for a similar statement. 
\end{rmk}

\begin{rmk}\label{rmk:cxc} {\rm  A similar statement than Theorem \ref{thm:main} holds for topologically cxc maps
with planar repellors, see \cite{kmp:ph:cxci} for
their definition and basic properties. 
But the complexity of the topology of the repellors and the lack of algebraic structure of such maps
require to develop other ingredients, 
so it will be explained elsewhere.}
\end{rmk}

\begin{corollary}[Bonk and Kleiner]\label{cor:main} A carpet group $G$ is virtually isomorphic to 
a convex-cocompact Kleinian group if and only if $\confdim_{AR} G< 2$.\end{corollary}

\Proof Carpets have essentially one embedding in the sphere up to homeomorphisms of the sphere.
Since boundary components do not separate the carpet, the group always has a planar
action. So Theorem \ref{thm:main} applies.\endp  

\begin{rmk}\label{rmk:bonkleiner} {\rm This corollary was announced
by Bonk and Kleiner in  2006 \cite{bonk:icm:qcgeom} where a sketch of the proof was given.
If the first step of our proof of Theorem \ref{thm:main} ---and hence of Corollary \ref{cor:main}--- 
is similar to theirs (filling-in the ``holes'' to reconstruct the Riemann sphere), the other
steps are different.
For Corollary \ref{cor:main}, Bonk and Kleiner exploit the specificity of carpets to extend the action to the whole sphere \cite{bonk:kleiner:merenkov}, whereas
here, we rely on  Hinkkanen-Markovic's characterization of M\"obius groups of the circle; this has the advantage to cover the cases
of other groups having a planar action.   Moreover, we rely on
the techniques developed in \cite{ph:sewing} to prove that
the carpet embeds quasisymmetrically in $\cbar$, see Corollary \ref{cor:qsembedcarpets}.
}
\end{rmk}

 
 Working a little more, we may obtain the following corollary from Theorem \ref{thm:main}:
 
 \REFCOR{cor:planaraction} Let $G$ be a torsion-free non-elementary hyperbolic group acting by homeomorphisms 
on $S^2$ as a convergence action. Let us assume that the restriction of its action to its limit set $\La_G (\ne S^2)$ is uniform.
If $\confdim_{AR}G<2$, then the action of $G$ is conjugate to that of a discrete group of M\"obius transformations.
\ENDCOR

Theorem \ref{thm:main2} and Corollary \ref{cor:main} provide us with the following equivalent statements in relation to 
 the Kapovich-Kleiner conjecture. 
 
 \begin{corollary}Let $G$ be a carpet group with a faithful and orientation preserving action on its boundary.
 The following are equivalent:
 \bit
 \item the group $G$ is isomorphic to a convex-cocompact Kleinian group;
 \item the group $G$ acts cellularly and geometrically on a CAT(0) cube complex;
 \item the Ahlfors regular conformal dimension of $G$ is strictly less than two.
 \eit
 \end{corollary}

\subsection{Outline of the paper and of the proofs}\label{sec:outlinepf} There are several key ingredients for the proofs of the main results of different
nature: analytical, topological and algebraic.

\subsubsection{Analytic aspects}
The proof of Theorem \ref{thm:main} is motivated by the following two theorems.

\begin{theorem}[Sullivan, \cite{DS3}]\label{thm:sullivan} A countable group of uniformly quasi-M\"obius transformations
on $\cbar$ is conjugate to a group of M\"obius transformations.\end{theorem}


\begin{theorem}[Bonk and Kleiner,  \cite{bonk:kleiner:qsparam}] \label{thm:bkinv} A metric $2$-sphere is quasisymmetrically equivalent to the Riemann sphere
if it is linearly locally connected and 2-Ahlfors-regular.\end{theorem}

A metric space $Z$ is {\em linearly locally connected}  if there is a constant $\la\ge 1$ such that, 
 for all $z\in Z$ and $R>0$,
\bit
\item[(LLC1)]  for all $x,y\in B(z,R)$ there is a continuum $E\subset B(z,\la R)$
which contains $\{x,y\}$;
\item[(LLC2)] for all  $x,y\notin B(z,R)$, there is a continuum $E\subset Z\setminus B(z,(1/\la) R)$
which contains $\{x,y\}$.
\eit

The main argument of the proof of Theorem \ref{thm:main}  consists in extending the given planar action into  
a uniform quasi-M\"obius action of the group on a metric sphere quasisymmetric to $\cbar$. 
This idea was explained in \cite{bonk:icm:qcgeom} for carpet groups. 
The construction of the sphere boils down to gluing infinitely many so-called quasidisks to the boundary
of the group as was suggested by Heinonen. The assumption on the Ahlfors-regular conformal dimension in all theorems of the paper 
is used here ---to build the sphere. If it is already known that the boundary embeds quasisymmetrically in the Riemann sphere, then this assumption
becomes superfluous.  It turns out that the basic step of gluing two  metric spaces together while controlling the geometry was done in \cite[Thm.\,1]{ph:sewing}.
The present construction generalizes this basic step in some sort of obvious but tedious way. 
The extension of the dynamics is also based on \cite{ph:sewing}. 

The next three sections are concerned with the analytical aspects of the proofs. 
In Section \ref{sec:sewing}, we make a systematic analysis
of gluing together countably many continua to a fixed continuum and 
study the properties which are inherited. 
The proofs are routine and detailed but have the advantage to be checkable. The main results
are Theorem \ref{prop:gsewing} and Theorem \ref{thm:geoprop}. These results are 
specialized in \S\,\ref{sec:planar} to planar sets. 
Section \ref{sec:ext} is a continuation of Section \ref{sec:sewing}: it is explained
how quasi-M\"obius maps can be extended when enlarging the space.

We have made the following choices for these sections: as mentioned above, we have decided to work in the very general setting
of gluing infinitely continua to a given one, even if we are primarily interested in gluing disks to a planar compact set along simple closed curves. There are  two main 
reasons for that. First, the construction is very general and there would be no genuine difference in the exposition: working directly with disks
would have essentially replaced Proposition \ref{thm:beurling-ahlfors} with the characterization of quasicircles of low-dimension from \cite{herron:meyer:qcircle} ---which
is not more elementary.
Purely topological statements have been 
separated from those which are metric in nature; the specific features  of 
planar sets only appear at the end in the uniformization process. 
Second, such constructions should have applications to other settings:   gluing disks along more complicated spaces than simple closed curves
(this is the case  for topological cxc maps, cf. Remark \ref{rmk:cxc}) and also to non-planar settings. 

Groups enter in Section \ref{sec:hyp}. We provide some background on convergence group actions, word hyperbolic groups
and their quasiconvex subgroups. We establish some properties for word hyperbolic groups which are known to hold for convex-cocompact Kleinian groups. 
In particular, we state ---without proof--- some specific properties of convex-cocompact Kleinian groups in the next Proposition \ref{prop:cck} 
which serve as a motivation for the present approach to their characterization.

\begin{proposition}\label{prop:cck} Let $G$ be a non-elementary convex-cocompact group. The following properties hold.
\ben
\item The limit set is porous (def. in \S\ref{sec:geoprop}) and Ahlfors regular of dimension $Q<2$.
\item Each connected component of $\La_G$ is linearly locally connected. 
\item For any $\de>0$, there are only finitely
many components of $\oo_G$  of diameter at least $\de$ (cf. Prop.\,\ref{prop:sierpgroup}).
\item The quotient space $\oo_G/G$ is a finite number of compact orbifolds (see Remark \ref{rmk:ahlforsfiniteness}).
\item Each non-trivial connected component of the boundary of a component of $\oo_G$ is a uniform quasicircle (cf. Prop.\,\ref{prop:sierpgroup}).
\een\end{proposition} 

These properties enable us to reconstruct the Riemann sphere from the boundary of a group by applying the results from \S\S\,\ref{sec:sewing}--\ref{sec:planar}.
After  exhibiting a way of extending convergence actions  to larger spaces in Theorems \ref{thm:extcva} and \ref{thm:extpr},  we conclude Section \ref{sec:hyp}  with the proof
of Theorem \ref{thm:main} and of Corollary \ref{cor:planaraction}. 

\subsubsection{Low dimensional topology}\label{sec:lowtop}

In general, there is no reason why a word hyperbolic group  with a planar boundary should admit a planar action  and it does not even need to
be isomorphic to a Kleinian group, see \cite[\S\,8]{kapovich:kleiner:lowdim}. Thurston's geometrization theory of $3$-manifolds will  circumvent
this difficulty, cf. \cite{thurston:bams, morgan:thurston, kapovich:book} .

Recall that a  Kleinian group  is a discrete subgroup of $\PP SL_2(\C)$ which we view as acting
both on hyperbolic $3$-space $\HH^3$ via orientation-preserving isometries and on the Riemann
sphere $\cbar$ via M\"obius transformations. 
As Poincar\'e observed, we may identify the Riemann sphere with the boundary at infinity of  the three-dimensional hyperbolic space \cite{poincare:memoire:kleineen}. Explicitly, let us consider   
the open unit ball in $\R^3$ as a model of  $\HH^3$ and
 the unit sphere $\SS^2$ for  the  Riemann sphere. One obtains
in this way an action of a Kleinian group $G$ on the closed unit ball. With this identification in mind, the group $G$ preserves the convex hull $\hbox{Hull}(\La_G)$ of its limit set in $\HH^3$.
The group  $G$ is  convex-cocompact if its action is cocompact on $\hbox{Hull}(\La_G)$.

When $G$ is torsion-free, we may associate a  $3$-manifold $M_G= (\HH^3\cup\oo_G)/G$, canonically 
endowed with a complete hyperbolic structure in its interior, which is called the {\it Kleinian manifold}.
The group $G$ is convex-cocompact  if and only if  $M_G$ is compact. Conversely, a compact $3$-manifold $M$ is
{\it hyperbolizable} if there exists a discrete subgroup of isometries $G$
such that $M$ is homeomorphic
to $M_G$ (this whole presentation rules out  tori in $\partial M$ since they are not relevant to the present work).
 We say that $M$ is {\it uniformized} by $G$. Note that $G$ is isomorphic to the fundamental 
group of $M$, and that it is necessarily word hyperbolic. Moreover, the boundary $\partial M$ is a union
of finitely many hyperbolic compact surfaces, cf. Proposition \ref{prop:cck}. When $M$ is orientable, then $G$ is a  convex-cocompact Kleinian group.

A compact {\it Haken three-manifold} is a manifold which can be constructed topologically by piecing together finitely many three-dimensional balls
along their boundaries. A compact manifold with non-empty boundary is always Haken. 
Thurston established  a uniformization theorem for Haken manifolds,  which we sate as follows.

\begin{theorem}[Thurston]\label{thm:thurston}
A 
compact irreducible Haken $3$-manifold with an infinite
word hyperbolic 
 fundamental group 
 is hyperbolizable.
\ENDTHM

If the orientable case is usually stated, see for instance  \cite[Thm 2.3]{thurston:bams} and \cite[Thm A']{morgan:thurston},
this is  not the case for non-orientable manifolds. 
It can be  deduced from the uniformization of orientable manifolds: by taking the orientable double cover,
we obtain a representation of its fundamental group as a group generated by a convex-cocompact Kleinian group
and an orientation-reversing quasiconformal involution; this group is thus uniformly quasiconformal,
hence conjugate to a group of M\"obius transformations according to Sullivan's straightening theorem, Theorem \ref{thm:sullivan}.

In particular, one obtains the following caracterization of convex-cocompact Kleinian groups among fundamental groups of $3$-manifolds \cite[Cor.\,4.9]{ctm:classification}.
\REFCOR{cor:thurston} A compact $3$-manifold  $M$ is homeomorphic to the Kleinian
manifold of a convex-cocompact Kleinian group $G$ with $\La_G\ne\cbar$ if and only if
$M$ is irreducible, orientable with non-empty boundary and its fundamental group is word hyperbolic.
\ENDCOR

The proofs of the main results of the paper will rely on the algebraic
structure of the group  to construct a finite index subgroup isomorphic to the fundamental group of a compact
Haken $3$-manifold. Theorems \ref{cor:main1} and  \ref{thm:main1} will then  follow from Thurston's uniformization theorem. 

\subsubsection{Planar topology}
The starting point of this approach i.e., the construction of a manifold, is Bowditch's JSJ decomposition,  which splits a word hyperbolic group into a graph of groups of three
different types: virtually cyclic groups, virtually free groups and so-called rigid groups --- which can be almost any
other word hyperbolic group with a planar boundary, see \S \ref{sec:jsj}. The limit sets of those groups are identified
from the structure of the local cut points of the boundary of the ambient group. In our setting of planar boundaries, a fine analysis of the topology 
of the rigid limit sets shows that the action of the corresponding stabilizer is the restriction of a well-defined planar action  (Prop.\,\ref{prop:typeIIIplanar}).
This is our next key ingredient, which is established using  planar topology.
It follows that each subgroup arising in the JSJ decomposition is, up to finite index, the fundamental group of a compact
3-manifold provided its Ahlfors-regular conformal dimension is strictly less than two. 
This provides us with a  characterization of Kleinian groups when the group admits a {\em regular JSJ decomposition} as defined therein.

\subsubsection{Algebraic aspects coming from geometric group theory}
Of course, it is not known if a word hyperbolic  group $G$ (with planar boundary) is always virtually torsion-free, and, even if this is the case, whether we can glue the pieces given by the JSJ decomposition 
together to build a manifold with fundamental group isomorphic to (a finite-index subgroup of) $G$.
Theorem \ref{thm:qcerf} shows that this last step follows from the separability of the quasiconvex subgroups of $G$ (abbreviated as the QCERF property), see \S\,\ref{sec:qcerf} for
the definition.
This property was first introduced in a topological context by Scott in order to desingularize immersed submanifolds into finite covers \cite{scott:almostgeo}.
This is the key ingredient for solving the so-called virtual Haken conjecture \cite{agol:virtualHaken}. In the present work, this property is used
to remove the obstructions preventing the construction of a manifold with fundamental group isomorphic to a finite index
subgroup of $G$. 

The QCERF property is known to hold for word hyperbolic groups acting on CAT(0) cube complexes thanks to
the works of Haglund and Wise \cite{haglund:wise:special} and of Agol \cite{agol:virtualHaken}. 
This is explained in Section \ref{sec:cccat} which is essentially descriptive: we recall the main results obtained by Wise, Agol and others on CAT(0) cube complexes and their role  concerning  the
accessibility of groups and the separability of subgroups. 
This is the last ingredient which provides us with the missing step for proving that a group with a planar boundary is virtually isomorphic to a convex-cocompact
Kleinian group under appropriate assumptions.




\subsubsection{Conclusion}
The proofs of the main results are then established in the last section. 
A consequence of the results established in the previous sections is that a group is virtually isomorphic to a convex-cocompact Kleinian group if the rigid vertices appearing in JSJ decompositions have their Ahlfors-regular conformal dimension strictly less than $2$
(Prop.\,\ref{prop:basechar}). This enables us to prove Theorems \ref{thm:main1} and  \ref{cor:main1}. 
Remark \ref{rmk:qivsqs} says that Theorems \ref{thm:mmain} and  \ref{thm:cmain}  are consequences of
Theorem \ref{cor:main1}. 


Assuming there are no $2$-torsion provides us with a finite hierarchy so that
the resulting subgroups cannot be split over elementary groups: the remaining non-trivial subgroups are carpet groups.
We may then prove 
Theorem  \ref{thm:main0} along  the same lines as above by induction on the length of 
the hierarchy.

By the previous results, the proof of Theorem \ref{thm:main2} reduces to the case of carpet groups
and groups with boundary homeomorphic to the sphere. So we may assume that we are given a convergence group
action on the sphere, uniform on its limit set. We first show that we may define for such a group an 
action on a CAT(0) cube complex
such that the stabilizers of hyperplanes are isomorphic to convex-cocompact Fuchsian groups. Splitting inductively
along those hyperplanes, one will obtain hyperbolic manifolds endowed with subsurfaces on their boundary which can
be glued together to prove that $G$ is virtually the fundamental group of a compact Haken manifold. The proof ends as above.

\subsection{Acknowledgements}
 I wish to express my gratitude to Michael Kapovich for Theorems \ref{thm:mmain}, \ref{thm:cmain} and \ref{cor:main1} and
for his guidance on quasi-isometric rigidity; after a first
version of this work circulated, he pointed it out to me that its contents  implied the quasi-isometric rigidity of
convex-cocompact Kleinian groups. I also feel particularly indebted to Cyril Lecuire who explained to me
many features of hyperbolic manifolds and with whom I have had many fruitful discussions, and to Fr\'ed\'eric Haglund
for his explanations on cubulated groups and for suggesting Theorem \ref{thm:main2}.
I am also very grateful to the following people for the many discussions I have had around
this project and for their encouragements: J.\,Aramayona, M.\,Boileau, M.\,Carrasco, T.\,Coulbois,
P.\,Derbez, V.\,Guirardel, A.\,Hilion, J.\,Los, L.\,Paoluzzi, J.P.\,Pr\'eaux, H.\,Short, 
N.\,Touikan, A.\,Yaman, and others I might have forgotten. 
This work was partially supported by the
ANR project ``GDSous/GSG'' no. 12-BS01-0003-01.


\section{Sewing infinitely many continua}\label{sec:sewing}
This section is devoted to spaces obtained by the procedure of gluing together infinitely many compact sets.
We first introduce the sewing procedure from a purely topological perspective, then we specialize this process
to metric spaces. In the last paragraph,  we study geometric properties of the new space inherited from each subcontinua.

There are similarities between this section and the next with the work of Merenkov and Wildrick \cite{merenkov:wildrick:koebe}, but with
some major differences concerning both the topology (local cutpoints) and the measure theory (how the mass is distributed), 
which make both contributions complementary rather than overlapping. 



\subsection{Topological sewing}
Let $X$ be a {\em continuum} i.e., a Hausdorff non-degenerate connected compact space. 
We assume that we are given a {\it null-sequence} i.e., an at most countable family
$\PPP$ of subcontinua with the following property: for any finite cover $\UUU$ of $X$,
for all but finitely many elements $K$ of $\PPP$, there exists $U\in\UUU$ with $K\subset U$.
We call $\PPP$ an {\em admissible collection of boundary components of $X$}.

For each $K\in\PPP$, we assume that we are given a continuum $L_K$ together with
an injective mapping $\psi_K:K\to L_K$.

To make this section more concrete, it can be read having in mind a degenerate carpet as defined in
Section \ref{sec:planar}  for $X$, the collection of Jordan curves appearing as the boundaries of the connected components of its complement in $S^2$ for 
$\PPP$, and closed Euclidean disks for the $L_K$'s  glued to $X$ along their boundary.

Set $$\Si= X \sqcup (\cup_{K\in\PPP} L_K)/\sim$$
where, for all $K\in\PPP$,  $z\in K$ is identified with $\psi_K(z)$;
note that a point $z$ may belong to  several boundary components. 
We define a topology on $\Si$ as follows: a basis of open sets of $\Si$ consists of those
sets $U$ such that 
\bit
\item[(T1)] $U\cap X$ is open in $X$, 
\item[(T2)] $U\cap L_K$ is open in $L_K$ for all $K\in\PPP$,
\item[(T3)] for all but finitely many components $K\in\PPP$ with $K\subset U$, one also has
$L_K\subset U$.\eit
One may check that these sets are stable under finite intersection, so we have indeed a basis
for a topology.

\REFPROP{prop:sicomp} With the notation above, the topological space
$\Si$ is Hausdorff and compact, and each embedding $X\hookrightarrow \Si$ and $L_K\hookrightarrow \Si$
is continuous. The connected components of $\Si\setminus X$ are in bijection with 
$\{$the connected components
of $L_K\setminus \psi_K(K)$, $K\in\PPP\}$.  \ENDPROP

\Proof 
The continuity of the embeddings follows from (T1) and (T2).

We now construct a collection of open neighborhoods for each
point of $\Si$.

Let us first consider $x\in X$. If $U_x$ is an open neighborhood of $x$ in $X$
and $W_x$ is the interior of a compact neighborhood of $\partial U_x$ disjoint from $x$, then the collection $\FFF_x$
of boundary components which intersect $\partial U_x$  and is not contained
in $W_x$ is finite since $\PPP$ is admissible. Let $U_x'$ be the complement in $U_x$ of the union
of boundary components $K\notin\FFF_x$ with $K\cap \partial U_x\ne\emty$; the set $U_x'$ is an open
neighborhood of $x$ in $X$: indeed, if $u\in U_x'$ and $N\subset U_x$ is a compact neighborhood of $u$, then
only finitely many components $K\in\PPP$ intersect both $\partial U_x$ and $N$, so their complement
in $N$ is a neighborhood of $u$ in $U_x'$.

For each $K\in\FFF_x$, $\psi_K(U_x\cap K)$ is  open in $\psi_K(K)$, so there exists
an open set $U_K\subset L_K$ such that $\psi_K(U_x\cap K)= U_K\cap \psi_K(K)$. 
For each $K\in\PPP\setminus\FFF_x$ such that $K\cap U_x'\ne\emty$, we let $U_K=L_K$.
For the other components, set $U_K=\emty$.
It follows that $$V_x= U_x'\cup (\cup_{K\in\PPP} U_K)$$
is an open neighborhood of $x$ in $\Si$.

If $x\notin X$, then there exists $K\in\PPP$ with $x\in L_K\setminus \psi_K(K)$. Let
$U_x$ be a neighborhood of $x$ in $L_K$, then $U_x\setminus \psi_K(K)$ is an open
neighborhood of $x$ in $\Si$. 

It follows easily that $\Si$ is Hausdorff.

 Let us now consider a covering of $\Si$ by open sets. We may as well
assume that each element satisfies (T1), (T2) and (T3).
Since $X$ is compact, we may extract a finite cover $\UUU_0$ of $X$.
The admissibility of $\PPP$ and condition (T3) imply that the set $Y=\Sigma\setminus \cup_{U\in\UUU_0} U$ intersects only finitely 
many elements $L_K$, each of which is compact by assumption. Therefore, one may extract a finite
cover for each of these sets and obtain a finite cover of $\Si$.

The last statement follows easily since the sets $L_K\setminus\psi_K(K)$, $K\in\PPP$, are pairwise disjoint
and their union forms $\Si\setminus X$.\endp

\subsection{Geometric sewing}


We prove a metric version of Proposition \ref{prop:sicomp}:

\REFTHM{prop:gsewing} Let $X$ be a metric continuum endowed with an admissible collection of
boundary components $\PPP$.
We assume the existence of $\De_0\ge 1$ and $\eta:\R_+\to\R_+$ such that, 
for each $K\in\PPP$, we are given a metric continuum $L_K$ and an $\eta$-quasisymmetric 
embedding  $\psi_K:K\to L_K$ such that $\diam\,L_K\le \De_0\diam\,\psi_K(K)$.

Then there exist a metric $d_{\Si}$ on $\Si$ compatible with its topology and a constant
$\De>0$ such that 
$(X,d_{\Si})$
is bi-Lipschitz to $X$, and,  for all $K\in\PPP$,  $(L_K,d_{\Si})$ is uniformly quasisymmetric
to $L_K$, $\diam_{\Sigma} L_K\le \De \diam_{\Sigma} K$ and there is a constant $c>0$ such that,
for all $y\in X$, $z\in L_K$, $z'\in L_{K'}$, $K\ne K'$, 
\begin{equation}\label{eq:flat0}
\begin{split}
& 
d_{\Sigma}(z,z') \ge c\inf \{d_{\Sigma}(z,x)+d_{\Sigma}(x,x')+d_{\Sigma}(x',z'),\ x\in L_K,x'\in L_{K'}\} 
\\
& 
d_{\Sigma}(z,y)  \ge c\inf \{d_{\Sigma}(z,x)+d_{\Sigma}(x,y),\ x\in L_K\}\,.
\end{split}\end{equation} \ENDTHM

The basic distortion bound for quasisymmetric maps is given by the following lemma
\cite[Prop.\,10.8]{heinonen:analysis}, which will be used throughout the paper:
\REFLEM{lma:qsbounds} Let $h:X\to Y$ be an $\eta$-quasisymmetric map
between compact metric spaces. For all $A,B \subset X$ with $A \subset B$ and $\diam\,B<\infty$,  we have $\diam\,h(B)<\infty$ and
\[ \frac{1}{2\eta\left(\frac{\diam\,B}{\diam\,A}\right)} \leq \frac{\diam  \,
h(A)}{\diam\,h(B)} \leq
\eta\left(2\frac{\diam\,A}{\diam\,B}\right)\,.\] \ENDLEM

For the proof, we will also use the following Ahlfors-Beurling type theorem \cite[Thm\,2]{ph:sewing}:

\REFPROP{thm:beurling-ahlfors} Let $(X,d_X)$ be a proper metric space
and $(Y,d_Y)$ a connected compact metric space. Let us assume that there
is an $\eta$-quasisymmetric embedding $f:Y\to X$ with $\diam_Y Y=\diam_X f(Y)$. 
Then  there is a metric $\hat{d}$ on $X$ such that
\ben
\item
$Id:(X,d_X)\to (X,\hat{d})$ is $\heta$-quasisymmetric;
\item  $Id:(X\setminus f(Y),d_X)\to (X\setminus f(Y),\hat{d})$ is locally
 quasisimilar: there is a finite constant $C \ge 1$ such that, 
for any $x\in X\setminus f(Y)$ and any $y,z\in B_X(x, d_X(x,f(Y))/2)$,
$$\frac{1}{C}\le \frac{\hat{d}(y,z)}{d_X(y,z)}\cdot\frac{d_X(x,f(Y))}{\hat{d}(x,f(Y))}\le C\,;$$
\item $f:(Y,d_Y)\to (X,\hat{d})$ is bi-Lipschitz onto its image: there exists $L\ge 1$ such that, 
for all $y_1,y_2\in Y$,
$$\frac{1}{L} d_Y(y_1,y_2)\le \hat{d}(f(y_1),f(y_2))\le d_Y(y_1,y_2)\,;$$
\item there is a constant $\De\ge 1$ such that
$$\frac{1}{\De}\diam (X,d_X)\le \diam \left(X,\hat{d}\right)\le \De \diam (X,d_X)\,.$$\een
All the constants involved and $\heta$ only depend on $\eta$.
\ENDPROP

The proof of Theorem \ref{prop:gsewing} consists in defining the metric from the necessary
conditions given by the conclusion of the statement and then to check it fulfills the 
requirements.

\Proof (Thm\,\ref{prop:gsewing}) For $K\in\PPP$, we rescale the metric on $L_K$ so that
$\diam\,K=\diam\,\psi_K(K)$. We apply Proposition \ref{thm:beurling-ahlfors} to $L_K (=X)$, $K(=Y)$ and $\psi_K:K\to L_K$.
Let $d_K (=\hat{d})$ be the metric
thus obtained. Note that the collection of maps $\{L_K\stackrel{Id}{\longrightarrow}(L_K,d_K)\}_{K\in\PPP}$
is uniformly quasisymmetric with distortion function $\heta$. Moreover, Lemma \ref{lma:qsbounds} implies
that,  for all $K\in\PPP$, 
\begin{equation}\label{eq:2} \diam (L_K,d_K)\le \De \diam (\psi_K(K),d_K)\,,\end{equation}
where $\De=2\heta(\De_0)$, and, there exists $L\ge 1$ such that, for all $K\in\PPP$
and for all $x,y\in K$, 
\begin{equation}\label{eq:1}
\frac{1}{L} d_X(x,y)\le d_K(\psi_K(x),\psi_K(y))\le d_X (x,y)\,.\end{equation}
In the sequel, we will omit $\psi_K$ when it leads to no confusion.

Let us define a quasimetric on $\Si$ as follows:
\bit
\item if $x,y\in X$, set $q(x,y)= d_X(x,y)$;
\item if $x,y\in L_K\setminus K$, for some $K\in\PPP$, set $q(x,y)=d_K(x,y)$;
\item if $x\in X$ and $y\in L_K\setminus K$ for some $K\in\PPP$, set $$q(x,y)=q(y,x)=\inf_{z\in K} \{d_X(x,z)+d_K(z,y)\}\,;$$
\item if $x\in L_{K_1}\setminus K_1$, $y\in L_{K_2}\setminus K_2$ for some $K_1\ne K_2\in\PPP$, set
$$q(x,y)=\inf_{(z_1,z_2)\in K_1\times K_2} \{ d_{K_1}(x,z_1)+d_X(z_1,z_2)+d_{K_2}(z_2,y)\}\,.$$
\eit

Set finally $$d_{\Si}(x,y)=\inf \sum_{i=0}^{N-1} q(x_i,x_{i+1})$$
over all finite chains $x_0,,\ldots x_N$ in $\Si$ with $x_0=x$, $x_N=y$.
We claim that $d_{\Si}$ is a metric comparable to $q$  which is compatible with the topology of $\Si$.

Let $x_0,\ldots, x_N$ be a chain. Inserting finitely many points if necessary in the chain
using the definition of $q$, we may assume that, for each $j\in\{0,\ldots, N-1\}$,
either there exists $K \in\PPP$ such that $q(x_{j},x_{j+1})=d_K(x_j,x_{j+1})$, or $q(x_j,x_{j+1})=d_X(x_j,x_{j+1})$.
Using the triangle inequality, we may assume that if $x_j\notin X$, $j<N$, then
$x_{j+1}\in X$; if furthermore $0<j<N$, then there is some $K\in\PPP$ with $x_{j\pm 1}\in K$ and so, with (\ref{eq:1}),
$$q(x_{j-1},x_j)+q(x_j,x_{j+1})\ge d_K(x_{j-1},x_{j+1})\ge \frac{1}{L} d_X(x_{j-1},x_{j+1})\,.$$
Therefore, one can extract a subchain $(y_j)_{0\le j\le M}$ with $y_0=x_0$ and $ y_M=x_N$ such that

\bit
\item[(a)] if $x_0,x_N\in X$, then
$$\sum_{j=0}^{N-1} q(x_j,x_{j+1}) \ge  \frac{1}{L}\sum_{j=0}^{M-1} d_X(y_j,y_{j+1})\\
 \ge \frac{1}{L}d_X(y_0,y_{M}) = \frac{1}{L} q(x_0,x_{N})\,;$$
\item[(b)] if $x_0\in L_K\setminus K$ for some $K\in\PPP$ and $x_N\in X$, then $y_1\in K$ and
\begin{eqnarray*} 
\sum_{j=0}^{N-1} q(x_j,x_{j+1}) & \ge & d_K(y_0,y_1) + \frac{1}{L}\sum_{j=1}^{M-1} d_X(y_j,y_{j+1}) \\
& \ge & d_K(y_0,y_1) + \frac{1}{L}d_X(y_1,y_{M}) \ge \frac{1}{L} q(x_0,x_N)\,.\end{eqnarray*}
\item[(c)] if $x_0,x_N\notin X$, then
\begin{eqnarray*} 
\sum_{j=0}^{N-1} q(x_j,x_{j+1}) & \ge & q(y_0,y_1) + \frac{1}{L}\sum_{j=1}^{M-2} d_X(y_j,y_{j+1}) + q(y_{M-1},y_M)\\
& \ge & q(y_0,y_1) + \frac{1}{L}d_X(y_1,y_{M-1}) + q(y_{M-1},y_M)\,.\end{eqnarray*}

If there is some $K\in\PPP$ such that $y_0,y_M\in L_K\setminus K$, then $y_1,y_{M-1}\in K$ and it follows from
(\ref{eq:1}) that
\begin{eqnarray*}  \sum_{j=0}^{N-1} q(x_j,x_{j+1}) & \ge & d_K(y_0,y_1) + \frac{1}{L}d_X(y_1,y_{M-1}) + d_K(y_{M-1},y_M)\\
& \ge & \frac{1}{L} d_K(y_0,y_M)= \frac{1}{L} q(x_0,x_N)\,.\end{eqnarray*}

If not, then
\begin{eqnarray*}  
\sum_{j=0}^{N-1} q(x_j,x_{j+1}) & \ge & \frac{1}{L}\left (q(y_0,y_1) + d_X(y_1,y_{M-1}) + q(y_{M-1},y_M)\right)\\
& \ge & \frac{1}{L} q(y_0,y_M)= \frac{1}{L} q(x_0,x_N)\,.\end{eqnarray*}
\eit

In either case, we have shown that 
\begin{equation}\label{eq:0}q(x_0,x_N) \ge d_{\Si}(x_0,x_N)\ge \frac{1}{L} q(x_0,x_N)\,.\end{equation}
This proves (\ref{eq:flat0}).
Since $L_K$ embeds in $\Si$ uniformly quasisymmetrically,
Lemma \ref{lma:qsbounds} implies the existence of $\De>0$ such that  $\diam_{\Sigma} L_K\le \De \diam_{\Sigma} K$ for
all $K\in\PPP$. 

Since $\PPP$ is a null-sequence and  (\ref{eq:2}) holds, it follows 
 that $d_{\Si}$ defines the topology of $\Si$.
 \endp

\begin{rmk}\label{distanceb} Note that if $x\in L_K\setminus \psi_K(K)$ for some $K\in\PPP$, 
then $d_{\Si}(x,X)$ is realized
by a point $y\in K$.\end{rmk}

\subsection{Geometric properties} \label{sec:geoprop}
We establish a series of properties of $\Si$ ---obtained by Theorem \ref{prop:gsewing}---
inherited from  the sets $X$ and $L_K$, $K\in\PPP$.  We focus on those properties which are known
to hold for the limit sets of convex-cocompact Kleinian groups and which are needed to apply Theorem \ref{thm:bkinv}, cf. Prop.\,\ref{prop:cck}.
The terms in the next statement will be defined below, when they are established.

\REFTHM{thm:geoprop} Under  the assumptions of Theorem \ref{prop:gsewing}, let
$(\Si,d_{\Si})$ be the given metric space. Then the following hold.
\ben
\item If $(X,d_X)$ and all $L_K$, $K\in\PPP$, satisfy the bounded turning property uniformly, 
then it also holds for $(\Si,d_{\Si})$. 
\item If $(X,d_X)$ and all $L_K$, $K\in\PPP$, are uniformly LLC, 
then $(\Si,d_{\Si})$ is LLC  quantitatively as well.
\item If $X$ is doubling and relatively doubling with respect to $\PPP$ and
if $\{L_K,\ K\in\PPP\}$ is uniformly doubling, then $(\Si,d_{\Si})$ is doubling quantitatively.

\item If $X$ is  relatively porous with respect to $\PPP$ and
if, for all $K\in\PPP$, $K$ is uniformly porous in $L_K$, then $X$ is porous in $(\Si,d_{\Si})$
quantitatively.
\item We assume that $X$ is doubling, doubling relative to $\PPP$ and porous relative to $\PPP$.
We also assume that every $K$ is uniformly porous in $L_K$. 
If  each $L_K$, $K\in\PPP$, is $Q$-Ahlfors-regular
with uniform constants for some $Q> 1$, and if 
$X$ is Ahlfors regular of dimension strictly less than $Q$,
then $(\Si,d_\Si)$ is $Q$-Ahlfors regular.
\een
\ENDTHM

We will first see how to extend local properties to global ones. Theorem \ref{thm:geoprop} will
follow at once.

Given a metric space $Z$ and a constant $c>0$, we say that $(X,\{L_\al\}_{\al\in A})$ 
is a {\it $c$-separating structure} of $Z$ if $X\subset Z$ is closed, $\{L_\al\}_{\al\in A}$ is a possibly infinite collection of closed
subsets of $Z$ which satisfies the following properties:
\bit
\item[(S1)] Setting $K_\al= L_\al \cap X$ and $\oo_\al =L_\al\setminus K_\al$ for $\al\in A$,  the collection $\{\oo_\al\}_{\al\in A}$ forms a 
partition of $Z\setminus X$ by open sets;
\item[(S2)] The following flatness condition holds: for all $y\in X$, $z\in \oo_\al$, $z'\in \oo_{\al'}$, $\al\ne \al'$,
\begin{equation}\label{eq:flat}
\left\{
\begin{split}
& 
d(z,z') \ge c\inf \{d(z,x)+d(x,x')+d(x',z'),\ x\in K_\al,x'\in K_{\al'}\} ;\\
& 
d(z,y)  \ge c\inf \{d(z,x)+d(x,y),\ x\in K_\al\} \,.
\end{split}\right.\end{equation}
\eit


Unless explicitly stated, we assume throughout this section 
that we are given a metric space $Z$ together with a separating
structure  $(X,\{L_{\al}\}_{\al\in A})$. 
The collection  $\PPP=\{K_\al,\ \al\in A\}$ denotes the boundary components of $X$. 
We also assume the existence of $\De>0$ such that $\diam\,L_\al \le \De \diam\, K_\al$ for
all $\al\in A$.

Note that if $z\in\oo_\al$ and $\de(z)=d(z,X)$, then 
\begin{equation}\label{eq:sep} B(z,c\de(z))\subset\oo_\al.\end{equation}

\subsubsection{Connectedness properties}
We establish connectedness properties which will prove that the space $\Si$ will be linearly locally connected, as required by  Theorem \ref{thm:bkinv}.

Recall that a metric space $Z$ satisfies the {\em bounded turning property} $\la$-(BT), for some $\la\ge 1$,
if any pair of points $x,y\in Z$
is contained in a continuum $E$ with $\diam\,E\le \la d_Z(x,y)$.

\REFPROP{prop:bt} Let $\la\ge 1$ be fixed.
If $X$ satifies the $\la$-(BT) condition and every subset $L_{\al}$  satisfies the $\la$-(BT) condition,
then $Z$ satisfies the $(\la/c)$-(BT) condition. 
\ENDPROP

\Proof Let $x,y\in Z$. If they both belong to either $X$ or to the same $L_\al$,  then  
the $\la$-(BT) property implies at once
the existence of a continuum $E$ containing $\{x,y\}$ with $\diam\, E \le \la d(x,y)$.

Let $x\in X$ and $y\in \oo_\al$ for some $\al\in A$.
By the separation condition (\ref{eq:flat}), 
there is some
$z\in L_\al$ such that $d(x,y)\ge c (d(x,z)+d(z,y))$.
There exist continua $E_1\subset X$ and $E_2\subset L_{\al}$ which contain $\{x,z\}$ and $\{z,y\}$ respectively and such that
$\diam\, E_1\le \la d(x,z)$ and $\diam\, E_2\le \la d(z,y)$. It follows that $E=E_1\cup E_2$ is a continuum containing
$x$ and $y$ and 
\begin{eqnarray*} \diam\, E & \le & \diam\, E_1 + \diam\, E_2\\
&\le & \la d(x,z) +\la  d(y,z)\\
& \le & ( \la/c) d(x,y)\,.\end{eqnarray*}

Assume $x\in \oo_{\al}$ and $y\in \oo_{\be}$ with $\al\ne \be$. 
Let $z_\al\in L_\al\cap X$ and $z_\be\in L_\be\cap X$ satisfy  $d(x,y)\ge c(d(x,z_\al)+d(z_\al,z_\be)+ d(z_\al,y))$.
There exist continua $E_0\subset X$, $E_\al\subset L_{\al}$ and $E_\be\subset L_{\be}$ which contain 
$\{z_\al,z_\be\}$, $\{x,z_\al\}$ and $\{z_\be,y\}$ respectively and such that
$\diam\, E_0\le \la d(z_\al,z_\be)$, $\diam\,E_\al\le \la d(x,z_\al)$ and $\diam\, E_\be\le \la d(z_\be,y)$. 
It follows that $E=E_0\cup E_\al\cup E_\be$ is a continuum containing $x$ and $y$ and 
\begin{eqnarray*} \diam\, E & \le &\diam\, E_0 + \diam\, E_\al + \diam\, E_\be\\
&\le & \la d(z_\al,z_\be) +\la d(x,z_\al) +\la  d(y,z_\be)\\
& \le & ( \la/c) d(x,y)\,.\end{eqnarray*}
\endp

Recall from the introduction the definition of  {\it linear local connectedness}  ($\la$-LLC for some $\la \ge 1$). 

\REFPROP{prop:LLC}  Let $\la\ge 1$ be fixed.
If $Z$ is $\la$-(LLC) at every point of $X$ and every
subset $L_\al$ satisfies the $\la$-(LLC) condition,
then $Z$ is LLC quantitatively as well.
\ENDPROP

\Proof By Proposition \ref{prop:bt}, the space $Z$ is $\la$-(BT) so (LLC1) holds. Let us now focus on (LLC2).
Let $z\in Z\setminus X$ and $r>0$. Let $\oo_\al$ be the component containing $z$.

If $d(z,X)>\frac{1}{2}\frac{r}{1+\la}$, then $B(z,\frac{1}{2}\frac{r}{\la(1+\la)})$ is disjoint from $X$ and the
(BT)-property implies that this ball is contained in $\oo_\al$, hence in $L_\al$. 
Therefore, if $z_1,z_2\notin B(z,r)$, then we may
find a continuum $K$ disjoint from $B\left(z,\frac{1}{2}\frac{r}{\la^2(1+\la)}\right)$ which joins $z_1$ and $z_2$.

If $d(z,X)\le\frac{1}{2}\frac{r}{1+\la}$, let $x\in X$ satisfy $d(z,X)=d(z,x)$. Note that
$$B\left(x,\left(1-\frac{1}{2}\frac{1}{1+\la}\right)r\right)\subset B(z,r)$$ and $$1-\frac{1}{2}\frac{1}{1+\la} = \frac{1 +2\la}{2(1+\la)}.$$
If $z_1,z_2\notin B(z,r)$, then we may
find a continuum $K$ disjoint from $B(x ,\frac{r}{\la}\frac{1+2\la}{2(1+\la)}))$ which joins $z_1$ and $z_2$.
Now, let us observe that
$$B\left(z,\left(\frac{1}{\la}\frac{1+2\la}{2(1+\la)}- \frac{1}{2}\frac{1}{1+\la}\right)r\right)
\subset B\left(x ,\frac{r}{\la}\frac{1+2\la}{2(1+\la)}\right)$$
so that $K$ is disjoint from $B(z,r/(2\la))$.\endp

\REFCOR{cor:LLC}   Let $\la\ge 1$ be fixed.
If $X$  and every
subset $L_\al$ satisfies the $\la$-(LLC) condition,
then $Z$ is LLC quantitatively as well.\ENDCOR

\Proof In order to apply Proposition \ref{prop:LLC}, it is enough to prove the following claim:
{\it there exists a constant $\la'\ge 1$ such that, 
for any $x\in X$, for any $r\in (0,\diam\, Z)$,  any point $z\in  (Z\setminus B(x,r))$ can be joined by a continuum
$E\subset  (Z\setminus B(x,r/\la'))$ to $X$. }

The claim and the LLC-property of $X$ imply that $Z$ is LLC at every point of $X$ quantitatively, so that
Proposition \ref{prop:LLC} applies. We now prove the claim.

Let  $x\in X$, $r>0$  and $z\notin B(x,r)$; there is some $\al \in A$ with $z\in\oo_\al$.
Set $$\kappa = \frac{1}{2(1+2\De \la)}\,.$$
If $d(x,K_\al )\ge  \kappa r$ then  $L_\al \cap B(x,c\kappa  r)=\emty$ according to (\ref{eq:flat}) and $L_\al$ is a continuum which
joins $z$ to a point of $K_\al\subset X$.
If $d(x,K_\al ) <  \kappa r$, let $x'\in K_\al$ be such that $d(x,x')= d(x,K_\al)$. 

It follows that $B(x',(1-\kappa) r)\subset B(x,r)$ so that $z\notin B(x',(1-\kappa) r)$.
Moreover,
$$\diam\,K_\al \ge \frac{1}{\De}\diam\, L_\al\ge  \frac{d(x',z)}{\De}\ge \frac{(1-\kappa)r}{\De}\,.$$
Therefore, we may find $y\in K_\al$ such that $d(x',y)\ge  \frac{(1-\kappa)r}{2\De}$.

Since $L_\al$ is LLC, there is a continuum $$E\subset \left(Z\setminus B\left(x', \frac{{(1-\kappa)r}}{2\la\De}\right)\right)$$
which joins $z$ to $y$.  Now, if $w\in E$, then $d(x,w)\ge d(w,x')- d(x',x)$
so that
$$d(x,w)\ge \frac{{(1-\kappa)r}}{2\la\De} - \kappa r \ge \frac{1}{4\la \De} r\,.$$
Therefore, the claim follows with $$\la'= \max\left\{ \frac{2(1+2\De \la)}{c},  4\la \De\right\}\,.$$
\endp

\subsubsection{Size properties}
For  Theorem \ref{thm:bkinv}, we also need to establish the Ahlfors-regularity of $\Si$, cf. Proposition \ref{prop:AR}.

A metric space $Z$ is {\em doubling} if there exists an integer $N$ such that
any set of finite diameter can be covered by at most $N$ sets of half its diameter.
This implies that, for all $\ep>0$, there exists $N_\ep$ such that any set $E$ of finite diameter
can be covered by $N_\ep$ sets of diameter bounded by $\ep\diam E$. 
We propose a relative notion of doubling:

\begin{defn}[Relative doubling condition] Let $X$ be a metric continuum with boundary components
$\PPP$. Then $X$ is {\em doubling relative to $\PPP$} if, for any $\ep>0$, there is some $N_\ep$ and
there exists $r_0>0$ such that, for any $x\in X$
and $r\in (0,r_0)$, there are at most $N_\ep$ components $K\in\PPP$ such that $B(x,r)\cap K\ne\emty$ and
$\diam (K\cap B(x,r))\ge \ep r$.\end{defn}

\REFPROP{prop:doubling} If $X$ is doubling and relatively doubling with respect to $\PPP$ and
if $\{L_\al,\ \al\in A\}$ is uniformly doubling, then $Z$ is doubling quantitatively.
\ENDPROP

\Proof 
Let us first consider $x\in X$, $r\in (0,\diam\,Z)$ and $\ep\in (0,1)$. We may cover $B(x,r)\cap X$ by
$M_1$ balls $\{B_j\}$ of radius $(\ep r/6)$, where $M_1$ only depends on the doubling condition of 
$(X,d_X)$ and on $\ep$. For each ball $B_j$, we add all the components $L_\al$, $\al\in A$, such that
$K_\al\cap B_j\ne\emty$ and $\diam\,K_\al\le (\ep r)/(6\De)$; we let $B_j'$ be the resulting set.
It follows that $\diam\,B_j'\le \ep r$. 

We are left with the components $K_\al\in\PPP$ with $K_\al\cap B(x,r)\ne\emty$ and $\diam\, K\ge (\ep r)/(6\De)$:
the relative doubling condition implies that there are at most $N$ such sets. Each of these sets can
be covered by $M_2$ sets of diameter at most $\ep r$  by the uniform doubling condition.

Therefore, for any $x\in X$, we may cover $B(x,r)$ by at most $M_1 + NM_2$ sets of diameter at most $\ep r$.

If $x\notin X$, then either $B(x,r)\subset L_\al$ for some $\al\in A$, and then we may use the doubling
condition of $L_\al$, or there is some $y\in X$ such that $B(x,r)\subset B(y,2r/c)$, cf. (\ref{eq:sep}).
Using $\ep=(c/4)$ above, we may cover $B(y,2r/c)$, hence $B(x,r)$ by a uniform number of sets of diameter
at most $(r/2)$.\endp

A subset $Y$ of a metric space $Z$ is said to be {\em porous} if there exists a constant $p>0$ such that
any ball centered at a point of $Y$ of radius $r\in (0,\diam Z]$ contains a ball of radius $p r$ disjoint from $Y$.
We propose a relative notion of porosity:

\begin{defn}[Relative porosity] Let $X$ be a metric continuum with boundary components
$\PPP$. Then $X$ is {\em porous relative to $\PPP$} if there exist a constant $p_X>0$ and a maximal size $r_0>0$ such that,
for any $x\in X$
and $r\in (0,r_0)$, there is at least one subcontinuum $K'$ of a boundary component $K\in\PPP$ such that
$K'\subset B_X(x,r)$, $K'\cap B_{X}(x,r/2)\ne\emty$ and $\diam_X K'\ge p_X r$. \end{defn}

\REFPROP{prop:porosity} If $X$ is  relatively porous with respect to $\PPP$ and
if, for all $\al\in A$, $K_\al$ is uniformly porous in $L_\al$, then $X$ is porous in $Z$
quantitatively.
\ENDPROP

\Proof  Let $x\in X$ and $r\in (0,r_0)$. 
There exist $\al\in A$ and $y\in K_\al$ such that $d(x,y)\le r/2$, and, if
$K'$ denotes the component of $K_\al\cap B(x,r)$ which contains $y$, then 
$\diam\, K'\ge p_X r$. 

We have $\diam\,L_\al \ge \diam\,K_\al \ge p_X r$ 
and $B(x,r)\supset B(y,p_X r /2)$.
It follows from the porosity of $K_\al$ in $L_\al$ that there is some $z\in \oo_\al$ such that
$B(z, cpp_X r/2)\subset \oo_\al\cap B(y,p_X r /2)$, so that $$B(z, c pp_X r/2)\subset B(x,r)\setminus X\,.$$
\endp


We may now establish the central result of this section.

\REFPROP{prop:AR}  Let $Q>0$. Let us assume that $Z$ is a doubling metric space and that $X$ is a porous Ahlfors regular 
compact subset of dimension strictly smaller than $Q$.
If every  subspace $L_\al$ is
 $Q$-Ahlfors regular  with uniform constants and $K_\al$ is uniformly porous in $L_\al$,
then $Z$ is $Q$-Ahlfors regular.\ENDPROP

The proof is exactly the same as \cite[Prop.\,2.18]{ph:sewing}. We repeat it for completeness.
Instead of asking $X$ to be Ahlfors-regular, it would have been enough to
require that its Assouad dimension be strictly less than $Q$, see \cite{luukkainen:assouad} for the definition and its main
properties.

\Proof
 Let us denote by $\mu$ the $Q$-Hausdorff measure in $Z$. 
If $x\notin X$, then there is some $\al\in A$ such that $x\in \oo_\al$;
we let $\de(x)= d(x,X)$. 


By assumption and (\ref{eq:sep}), we have $\mu(B(x, c r\de(x)))\asymp (r\de(x))^Q$ for all $r\in (0,1)$.


Let us fix a point $x\in X$ and $r>0$. 
Since $X$ is porous in $Z$, a constant $p>0$ exists such that $B(x,r)$ 
contains a ball $B(y,pr)$ 
disjoint from $X$. Therefore $pr\le \de(y)$ so 
$$ \mu(B(y,pr)) \ge \mu(B(y.c pr) )\gtrsim  r^Q  \quad \hbox{and}  \quad \mu(B(x,r))\ge \mu(B(y,pr)) \gtrsim r^Q\,.$$

For the converse inequality, we first note that since  the dimension of $X$ is strictly less than $Q$,
we have $\mu(X)=0$. Thus it is enough to bound $\mu(B(x,r)\setminus X)$.
We cover $B(x,r)\setminus X$ by balls $B(z, c\de(z)/10)$. We extract an at most countable subfamily 
$B(z_j,c\de_j/10)$ of pairwise disjoint
balls such that $B(x,r)\setminus X\subset \cup B(z_j, c\de_j/2)$ \cite[Thm\,1.2]{heinonen:analysis}.


Denote by $A_n$ the set of centers $(z_j)$ such that $r e^{-(n+1)}< \de(z_j)\le re^{-n}$. It follows that if
$z_j\in A_n$, then $\mu(B(z_j,c\de_j/2))\asymp \de(z_j)^Q\asymp r^Q e^{-Qn}$.

For each $z_j$, choose a point $x_j\in X$ such that $\de(z_j)= d(z_j,x_j)$. Since $Z$ is doubling 
 and
the balls $\{B(z_j,c\de_j/10)\}_j$ are disjoint, 
the nerve of the family of balls $\{B(x_j,\de(z_j)),\ z_j\in A_n\}$  has uniformly bounded valence $V$ 
(independent from $n$). 
Therefore,
we may split this family of balls into $V+1$ families of pairwise disjoint balls.
Since $\dim X<Q$,
the number of balls involved in $A_n$ is bounded by 
$e^{Q n}\te^n$, up to a factor
(which depends on $V$), for some $\te\in (0,1)$. 
Thus
$$\sum_{A_n} \mu(B(z_j,c\de_j/2))\lesssim e^{Q n}\te^n \left(e^{-n}r\right)^Q\lesssim \te^n r^Q.$$

Therefore  $$\mu(B(y,r))\le \sum_{n\ge 0} \sum_{A_n} \mu(B(z_j,\de_j/2))\lesssim \sum_{n\ge 0}\te^n r^Q\lesssim r^Q.$$

\bigskip

Let us consider a point $z\in Z\setminus X$, and let $x\in X$ be such that $\de (z)=d(z,x)$. 
If $r\in [c\de(z) ,2\de(z)]$, then 
$$\mu(B(z,r))\ge \mu(B(z,c\de(z)))\gtrsim \de(z)^Q\gtrsim r^Q.$$
On the other hand, 
$$\mu(B(z,r))\le \mu(B(x,2r))\lesssim r^Q.$$

\medskip

If $r\ge 2\de(z)$, then $B(x,r-\de(z))\subset B(z,r)\subset B(x,r+\de(z))$ with $r+\de(z)\le (3/2)r$ and $r-\de(z)\ge r/2$, so 
$\mu(B(z,r))\asymp r^Q$.\endp

\subsubsection{Geometric properties of the sewn space}
We prove Theorem \ref{thm:geoprop}. Theorem \ref{prop:gsewing} implies that $(X,\{L_K,\ K\in\PPP\})$
is a separating structure for $\Si$ and that there is  some  $\De>0$ such that $\diam_\Si L_K \le \De \diam_\Si K$ for
all $K\in\PPP$.

\noindent{\bf 1.} 
It follows from Lemma \ref{lma:qsbounds}
that we may assume that $(X,d_X)$ and all $(L_K,d_K)$, $K\in\PPP$, satisfy $\la$-(BT) for a fixed $\la\ge 1$.
Therefore, Proposition \ref{prop:bt} applies and $\Si$ has bounded turning.\endp

\noindent{\bf 2.} 
We recall that the LLC property is preserved quantitatively under quasisymmetric mappings 
\cite[Chap.\,15]{heinonen:analysis},
so that we may assume that $\{(L_K,d_\Si),\ K\in\PPP\}$ and $(X,d_\Si)$ are $\la$-LLC. 
Therefore, Corollary \ref{cor:LLC} applies and $\Si$ is LLC as well.\endp

\noindent{\bf 3.} We recall that the doubling condition is preserved quantitatively under quasisymmetric mappings
\cite[Chap.\,15]{heinonen:analysis} 
so that $\{(L_K,d_\Si),\ K\in\PPP\}$ is uniformly doubling. 
Since $(X,d_X)$ is bi-Lipschitz to $(X,d_\Si)$, it is also doubling and relatively doubling with respect to $\PPP$.
Proposition \ref{prop:doubling} applies and we may conclude that $\Si$ is doubling.\endp

\noindent{\bf 4.} We recall that the porosity of a subset is preserved quantitatively under quasisymmetric mappings
so that the sets $K$, $K\in\PPP$, are uniformly porous in $(L_K,d_\Si)$ and 
$(X,d_\Si)$ is also relatively porous to $\PPP$. 
Therefore we may apply Proposition \ref{prop:porosity} and conclude that $X$ is porous in $\Si$.\endp

\noindent{\bf 5.} 
We know from above that the space $(\Si,d_\Si)$ is doubling and that $X$ is also porous in $\Si$. 
Since $(X,d_X)$
is bi-Lipschitz to $(X,d_\Si)$, we get that $(X,d_{\Si})$ is also Ahlfors regular of  dimension strictly less than $Q$.
Finally, \cite[Prop.\,2.18]{ph:sewing} implies that $(L_K,d_\Si)$, $K\in\PPP$, are uniformly $Q$-Ahlfors regular. 
Therefore, Proposition \ref{prop:AR} applies.\endp

 \section{Planar continua}\label{sec:planar}

Let $X$ be a planar locally connected continuum, and let $\vp:X\to S^2$ be a topological
embedding. Note that there may be embeddings which are not compatible in the sense that if
$\vp_1,\vp_2:X\to S^2$ are two embeddings, then $\vp_2\circ\vp_1^{-1}$ might not be the restriction
of a selfhomeomorphism of the sphere.

We let $\PPP$ denote those subcontinua $K\subset X$ such that
$\vp(K)$ bounds a connected component of $S^2\setminus \vp(X)$.

The main result of this section is the following:

\REFTHM{thm:qsembed} Let $X$ be a locally connected planar metric continuum with boundary components $\PPP$
provided by an embedding $\vp:X\to S^2$.
We assume that $X$ is LLC, Ahlfors regular of dimension $Q<2$, relatively doubling and porous with
respect to $\PPP$, and that we are given uniformly quasisymmetric gluing maps $\psi_K:K\to L_K$ for all $K\in\PPP$,
where the sets $L_K$ are continua such that (a) $L_K\setminus \psi_K(K)$ is homeomorphic to an open disk; (b) all the sets
$L_K$ are uniformly Ahlfors regular of dimension $2$;  (c) $\psi_K(K)$ is uniformly porous in $L_K$; and (d) $\diam L_K\le
\De_0\diam\,\psi_K(K)$ for some universal constant $\De_0\ge 1$.

Then the space $\Si$ given by Proposition \ref{prop:sicomp} and Theorem \ref{prop:gsewing} is quasisymmetric
to $\cbar$ and there exists a quasisymmetric embedding $\phi:X\to \cbar$ compatible with $\vp$.\ENDTHM 

\subsection{Topological uniformization}
 Let $X$ be a locally connected planar continuum with boundary components $\PPP$
provided by an embedding $\vp:X\to S^2$. For each $K\in\PPP$, we assume that we are given
topological embeddings $\psi_K:K\to L_K$ for all $K\in\PPP$,
where the sets $L_K$ are continua such that  $L_K\setminus \psi_K(K)$ is homeomorphic to an open disk.

For each $K\in\PPP$, we denote by $U_K$ the component of $S^2\setminus \vp(X)$
with $\partial U_K=\vp(K)$. It follows from the Torhorst theorem \cite[Thm.\,VI.2.2]{whyburn:analytic_topology} 
that $K$ is itself locally connected. Since $\vp(X)$ is connected, it follows that $U_K$ is homeomorphic
to an open disk so there exists a homeomorphism $\vp_K:L_K\to\overline{U_K}$ such that 
$\psi_K|_K= (\vp_K^{-1}\circ \vp)|_K$.

\REFPROP{prop:homeo} The topological space $\Si$ defined as above is homeomorphic to $S^2$.\ENDPROP

\Proof
We first note that since $\vp(X)$ is a locally connected continuum of the sphere,  \cite[Thm.\,VI.4.4]{whyburn:analytic_topology} 
implies that $\PPP$ is an admissible collection of boundary components. Hence $\Si$ is compact by Proposition \ref{prop:sicomp}. 

Define $\phi:\Si\to S^2$ as follows:
on $ X$, set $\phi=\vp$; on $L_K$, $K\in \PPP$,
set $\phi= \vp_K$. 
Note that, if  $z\in K$ for some $K\in\PPP$, then
$$\phi(z)=\vp(z)= \vp_K\circ (\vp_K^{-1}\circ\vp)(z)=\phi ( \psi_K(z))$$
so that $\phi$ is well-defined.

Let us prove that $\phi$ is a homeomorphism.
First, $\phi$ is a bijection by construction.
Let us now prove that $\phi$ is continuous: let $U$ be an open subset of $S^2$ and let
$V=\phi^{-1}(U)$. On the one hand, one has 
$$V\cap X=\phi^{-1}(U)\cap \phi^{-1}(\vp(X))=\phi^{-1}(U\cap\vp(X))=\vp^{-1}(U)$$
so that $V\cap X$ is open, hence (T1) is true. On the other hand, let $K\in\PPP$, then
$$V\cap L_K= \phi^{-1}(U)\cap \vp_K^{-1}(\overline{U_K})=\vp_K^{-1}(U\cap \overline{U_K})$$
so $V\cap L_K$ is also open, so (T2) holds. By \cite[Thm\,VI.4.4]{whyburn:analytic_topology}, 
for all $\de >0$, there are only finitely
many components of $S^2\setminus\vp(X)$ with diameter at least $\de$; therefore,
we may assume that only finitely many components of $S^2\setminus \vp(X)$ intersect $\partial U$,
since any open set can be described as an at most countable union of such open sets.
Hence $V$ satisfies (T3) as well.
Therefore $\phi$ is continuous.
Since $\Si$ is compact, this implies that $\phi$ is a homeomorphism.\endp

\subsection{Quasisymmetric embedding} The embedding will be obtained thanks to Bonk and Kleiner's uniformization Theorem \ref{thm:bkinv}.


\Proof (Theorem  \ref{thm:qsembed})
By the assumption (a), Proposition \ref{prop:homeo} implies that $\Si$ is homeomorphic
to $S^2$. The assumption (d) implies that we may endow $\Si$ with a metric $d_\Si$ which enjoys the properties
given by Theorem \ref{prop:gsewing}. 
 The assumption (c), 
the LLC-assumptions  and the Ahlfors-regularity assumptions enable us to apply  Theorem  \ref{thm:geoprop}
to conclude  that $(\Si,d_\Si)$ is an LLC $2$-Ahlfors regular sphere.

It now follows from Theorem \ref{thm:bkinv}  
that $\Si$ is quasisymmetrically equivalent to $\cbar$,
so $X$ embeds quasisymmetrically into $\cbar$. This ends the proof of Theorem \ref{thm:qsembed}.\endp

\begin{rmk} {\em If the space $X$ is not assumed to be LLC but only to satisfy a weaker connectedness  condition, 
then the linear local connectedness of the sphere $\Si$, hence the existence of a quasisymmetric embedding of $X$ into $\cbar$,
might depend on the choice of the boundary components. A finer analysis would thus be required. } \end{rmk}

\subsection{Carpets} Define a {\em carpet} as 
a planar, one-dimensional, connected, locally connected compact space with no local cut point;
any such space is homeomorphic to the Sierpi\'nski carpet
and admits a unique embedding up to postcomposition by a selfhomeomorphism of the sphere \cite{whyburn:sierpinski}. 
It follows that the collection of boundary components is canonically defined, and that they are pairwise
disjoint simple closed curves. We call them {\em peripheral circles}, and we say that they are {\em uniform quasicircles}
if they are the images of the unit circle by a quasisymmetry under a uniform distortion function. 
By extension, a {\it degenerate carpet} will be  a  one-dimensional, connected 
locally connected compact space
homeomorphic to the complement of a union of disjoint open disks (their closures may intersect). 
In this case,  there may be several non-equivalent embeddings in the sphere.

\REFLEM{lma:degcarpet} Let $X$ be a one-dimensional, connected 
locally connected planar compact space with no global cut point. Then any embedding of  $X$ in $S^2$ is  a
degenerate carpet.\ENDLEM

\Proof  Let us assume that $X$ is already embedded into $S^2$, that we identify with the Riemann sphere. 
Since it is one-dimensional, it has no interior in $S^2$. 
We now prove that the boundary of any component of $\oo$ of $S^2\setminus X$  is a  Jordan curve. This will establish that $X$
is a degenerate carpet.
We consider a conformal map $h:\DD\to\oo$. Since $\partial\oo$ is contained in a locally connected 
compact set (disjoint from $\oo$), Carath\'eodory's theorem implies that $h$ admits
a continuous and surjective extension $h:\overline{\DD}\to \overline{\oo}$. 
If $\oo$ is not a Jordan domain, then there are two rays in $\DD$ which
are mapped to a Jordan curve in $\overline{\oo}$ which separates $\partial\oo$, hence $X$.
But $X$ has no (global) cut point. Therefore, $\oo$ is a Jordan
domain and $X$ is a degenerate carpet. \endp

\REFCOR{cor:qsembedcarpets} Let $X$ be a metric degenerate carpet with $\confdim_{AR} X <2$  
endowed with boundary components which are assumed to be uniform quasicircles. 
We assume that $X$ is LLC,  relatively doubling and porous with
respect to the boundary components.
Then $\Si$ is quasisymmetric to $\cbar$ and 
there exists a quasisymmetric embedding $\phi:X\to \cbar$ compatible with the boundary components.
\ENDCOR  

\Proof We may choose an Ahlfors regular metric in the gauge of $X$ of dimension $Q<2$. 
For each $K\in\PPP$, there exists a uniform quasisymmetric homeomorphism $\psi_K:K\to\SS^1\,(\subset\overline{\DD})$.
Note that $\overline{\DD}$ is $2$-Ahlfors regular and LLC, that $\SS^1$ is porous in $\overline{\DD}$ and that 
$\diam\,\DD\le \diam\,\SS^1$.
Therefore, Theorem \ref{thm:qsembed} applies.\endp


\section{Extension of maps}\label{sec:ext}
In this section, we show how homeomorphisms between the different
sets can be glued together to yield a global homemorphism of $\Si$.

Let $X$ be a continuum endowed with  an 
admissible collection of boundary components $\PPP$.
For each $K\in\PPP$, we assume that we are given a continuum $L_K$ 
together with an injective mapping $\psi_K:K\to L_K$; and we consider
as above $$\Si= X \sqcup (\cup_{K\in\PPP} L_K)/\sim$$

\subsection{Global homeomorphisms}
The starting point is a collection of homeomorphisms $(h_K)_{K\in\PPP}$ and a homeomorphism
$h_X:X\to X$ such that, for all $K\in\PPP$, $h_X(K), h_X^{-1}(K)\in\PPP$. We assume that the following compatibility condition holds:
for all $K\in \PPP$, 
$h_K(\psi_K(K))=\psi_{h_X(K)}(h_X(K))$ and 
$$h_X|_K= (\psi_{h_X(K)}^{-1}\circ h_K\circ \psi_K)|_K\,.$$


\REFLEM{lma:homeoglobal} The map $h:\Si\to\Si$ defined by
$h(x)=h_X(x)$ if $x\in X$ and $h(x)=h_K(x)$ if $x\in L_K$
is a well-defined homeomorphism.\ENDLEM

\Proof We let the reader check that the compatibility condition implies that 
$h:\Si\to\Si$ is a well-defined bijection. 

Let $U\subset \Si$ be an open set which satisfies (T1), (T2) and (T3). 
Then $U\cap X$ is open in $X$ 
so $h^{-1}(U\cap X)=h_X^{-1}(U\cap X)$ is open in $X$ since $h_X$ is continuous.
Similarly, for any $K\in\PPP$, $U\cap L_K$ is open in $L_K$ 
so $h^{-1}(U\cap L_K)=h_{h_X^{-1}(K)}^{-1}(U\cap L_K)$ is open in 
$L_{h_X^{-1}(K)}$ since $h_{h_X^{-1}(K)}$ is continuous.

We note that $K\subset U$ if and only if $h_X^{-1}(K)\subset h_X^{-1}(U\cap X)$,
and $L_K\subset U$ if and only if $h^{-1}(U\cap L_K)\subset h^{-1}(U)$.
So, if for all but finitely many components $K\in\PPP$ with the property
that $K\subset U$ one has $L_K\subset U$, then the same is true for boundary
components $K\in\PPP$ with $K\subset h^{-1}(U)$. So $h^{-1}(U)$ is open,
and $h$ is continuous.

Since $h$ is also a bijection and $\Si$ is compact according to Proposition \ref{prop:sicomp}, it follows
that $h$ is a homeomorphism as well.\endp

\subsection{Quasi-M\"obius versus quasisymmetric maps}
Quasi-M\"obius  and quasisymmetric maps were defined in the introduction. 
We record the following relationships, 
see \cite{vam} for a proof.
\begin{proposition}\label{prop:qmvsqs} Let $f:Z\to Z'$ be a homeomorphism between proper metric spaces.
\begin{itemize}
\item[(i)] 
If $f$ is $\eta$-quasisymmetric then $f$ is also $\te$-quasi-M\"obius, where $\te$
only depends on $\eta$.

\item[(ii)]
If $f$ is $\te$-quasi-M\"obius, then $f$ is locally $\eta$-quasisymmetric, where $\eta$
only depends on $\te$.

\item[(iii)]
Let us assume that $f$ is $\te$-quasi-M\"obius. If $Z$ and $Z'$ are unbounded, then
$f$ is $\te$-quasisymmetric. If $Z$ and $Z'$ are compact, then assume that
there are three points $z_1,z_2,z_3\in Z$, such that $|z_i-z_j|\ge \diam\,Z/\la$ and
 $|f(z_i)-f(z_j)|\ge \diam\,Z'/\la$ for some $\la>0$, then $f$ is
$\eta$-quasisymmetric, where $\eta$ only depends on $\te$ and $\la$.
\end{itemize}
\end{proposition}

\REFCOR{cor:qmvsqs} Let $f:Z\to Z'$ be a $\te$-quasi-M\"obius homeomorphism between
proper metric spaces. If there is a closed and connected subset $Y\subset Z$ with at least three points
such that $\diam\,Y\ge \De \diam\,Z$, $\diam\,f(Y)\ge \De\diam\,Z'$ 
and such that $f|_Y$ is $\eta$-quasisymmetric, 
then $f$ is $\heta$-quasisymmetric on $Z$ where
$\heta$ only depends on $\De$, $\eta$ and $\te$.
\ENDCOR

\Proof If $Z$ is unbounded, then there is nothing to prove. If not,
we may pick three points $z_1,z_2,z_3\in Y$, such that $|z_i-z_j|\ge \diam\,Y/3$;
by Lemma \ref{lma:qsbounds}, we also have $|f(z_i)-f(z_j)|\ge \diam\,f(Y)/\la$ for some
constant $\la$ which only depends on $\eta$. But the assumption on the embeddings 
$Y\hookrightarrow Z$  and $f(Y)\hookrightarrow Z'$ implies that  $|z_i-z_j|\ge \De\diam\,Z/3$
and  $|f(z_i)-f(z_j)|\ge \De\diam\,Z'/\la$ for all $i\ne j$. Therefore, Proposition \ref{prop:qmvsqs}
applies.\endp

\subsection{Gluing quasi-M\"obius maps together}
We assume that $X$ is now a metric continuum equipped with an admissible collection of
boundary components $\PPP$ and we suppose
the existence of $\De_0\ge 1$ and $\eta:\R_+\to\R_+$ such that, 
for each $K\in\PPP$, we are given a metric continuum $L_K$ and an $\eta$-quasisymmetric 
homeomorphism  $\psi_K:K\to L_K$ such that $\diam\,L_K\le \De_0\diam\,\psi_K(K)$.
We let $d_{\Si}$ be the metric given by Theorem \ref{prop:gsewing} on $\Si$.
We keep the same notation as in Lemma \ref{lma:homeoglobal}.

In this section, we prove

\REFTHM{prop:gluqm} Given a \homeo $\te:\R_+\to\R_+$, there
exists another \homeo $\te':\R_+\to\R_+$ such that, if 
each $h_K$ and $h$ are $\te$-quasi-M\"obius maps, then 
$h:(\Si,d_{\Si})\to (\Si,d_{\Si})$ is $\te'$-quasi-M\"obius.\ENDTHM

It is convenient to see $\Si$ as an unbounded space in order
to transform quasi-M\"obius maps into quasisymmetric ones. 
The following lemma makes the job. It may be considered as a converse
construction of Bonk and Kleiner given in \cite[Lma\,2.2]{bonk:kleiner:rigidity}.
We use the same arguments for the proof.

\begin{lemma}\label{lma:stereo}
Let  $(X,d,w)$ be a marked complete metric space.
There exists a complete metric $\hat{d}$ on $X\setminus\{w\}$
such that  $Id:(X\setminus\{w\},d)\to (X\setminus\{w\},\hat{d})$ is $\te$-quasi-M\"obius 
with $\te(t)= 16t$.
and $$ \frac{1}{4}  \frac{d(x,y)}{\de(x)\de(y)} \le \hat{d}(x,y)\le \frac{d(x,y)}{\de(x)\de(y)}$$ where $\de(x)= d(x,w)$.
\end{lemma}

\Proof Set, for $x,y\in X\setminus\{w\}$, $$q(x,y)=  \frac{d(x,y)}{\de(x)\de(y)}$$
and let  $$\hat{d}(x,y)= \inf \sum_{i=0}^{N-1} q(x_i,x_{i+1})$$
over all finite chains $x_0,\ldots, x_N$ in $X\setminus\{w\}$ with $x_0=x$, $x_N=y$.
By definition, $\hat{d}(x,y) \le q(x,y)$ holds. 

Without loss of generality, we may assume that $\de(x)\le \de (y)$; let us  fix
a chain $x_0,,\ldots x_N$  with $x_0=x$, $x_N=y$. We consider two cases.

If $\de(x_j)\le 2\de(x)$ for all $j\in\{0,\ldots, N\}$, then 
\begin{eqnarray*}
 \sum_{i=0}^{N-1} q(x_i,x_{i+1}) & \ge & \frac{1}{4\de(x)^2}\sum_{i=0}^N d(x_i,x_{i+1})\\
 &\ge & \frac{1}{4\de(x)^2} d(x,y)\ge \frac{1}{4}  q(x,y)\,.\end{eqnarray*}
 
We now assume that there is some $j>0$ with $\de(x_j)\ge 2\de(x)$.
Let us observe that $$q(x,y) \le \frac{\de(x) + \de(y)}{\de(x)\de(y)} \le \frac{2}{\de(x)}$$ and that,
for $u,v\in X\setminus\{w\}$,
 $$\left| \frac{1}{\de(u)}-\frac{1}{\de(v)} \right|= \frac{|\de(v)-\de(u)|}{\de(u)\de(v)}\le q(u,v)\,.$$
 Therefore, 
 \begin{eqnarray*} \sum_{i=0}^{N-1} q(x_i,x_{i+1}) & \ge &  \sum_{i=0}^{j-1}\left| \frac{1}{\de(x_{i})}-\frac{1}{\de(x_{i+1})} \right|\\
& \ge & \left| \frac{1}{\de(x)}-\frac{1}{\de(x_{j})} \right|\\
&\ge & \frac{1}{2 \de(x)} \ge \frac{1}{4}q(x,y)\,.\end{eqnarray*}
This proves that $\hat{d}$ is indeed a complete metric. 
The fact that the identity is quasi-M\"obius follows at once from the formulae. \endp 

We will need to control the relative diameters of the sets $L_K$ and $K$ in this new metric:

\begin{lemma}\label{lma:diamstereo}
Let $X$ be a complete metric space, $w\in X$ and $K\subset L\subset X\setminus \{w\}$ be such that
$\diam\,L\le \De \diam\,K$ and there is a constant $c>0$ such that, for all $x\in L$,
$$ d(x,w) \ge c \inf\{ d(x,y)+ d(y,w),\ y\in K\}\,;$$ then there is some constant $\hat{\De} >0$ which
only depends on $\De$ and $c$ such that $$\diam \,(L,\hat{d})\le \hat{\De}\diam\,(K,\hat{d})$$
where $\hat{d}$ is the metric on $X\setminus\{w\}$ given by Lemma \ref{lma:stereo}.\end{lemma}

\Proof We set $\de(x)=d(x,w)$ and $q(x,y)=d(x,y)/(\de(x)\de(y))$ as above, and let $\de(K)=d(w,K)$.

Let us first estimate $\diam_q K=\sup\{q(x,y),\ x,y\in K\}$; 
pick $x\in K$ such that $\de(x)=\de(K)$ and let $y\in K$ be
such that $d(x,y)\ge (1/2)\diam\,K$. It follows that $\de(y)\le \diam\,K+ \de(K)$ so that
\begin{eqnarray*}
 \diam_q K \ge q(x,y) & \ge  & \frac{1}{2} \frac{\diam\,K}{\de(K)\cdot (\de(K)+\diam\, K)}\\ 
&\ge & \frac{1}{4} \frac{\diam\,K}{\de(K)\cdot\max \{ \de(K), \diam\, K\}}\\
& \ge & \frac{1}{4} \min\left\{\frac{\diam\,K}{\de(K)^2}, \frac{1}{ \de(K)}\right\}\,.\end{eqnarray*}

We note that the assumptions imply that, for all $x\in L$,  $\de(x)\ge c\cdot\de(K)$. 

On the one hand, we have
$$\diam_q L\le \frac{\diam\,L}{c^2\de(K)^2} \le \frac{\De}{c^2} \frac{\diam\,K}{\de(K)^2}\,.$$
On the other hand, if $x,y\in L$, then
$d(x,y)\le \de(x)+\de(y)\le 2\max\{\de(x),\de(y)\}$
so that $$q(x,y) \le \frac{2}{\min\{\de(x),\de(y)\}} \le \frac{2}{c\de(K)}$$
hence
\begin{eqnarray*}
\diam_q L & \le &\min\left\{\frac{\De}{c^2} \frac{\diam\,K}{\de(K)^2},  \frac{2}{c\de(K)}\right\}\\
& \le & \frac{2\De}{c^2} \min \left\{\frac{\diam\,K}{\de(K)^2},  \frac{1}{\de(K)}\right\}\\
& \le & \frac{8\De}{c^2} \diam_q K\,.\end{eqnarray*}
Therefore, $$\diam \,(L,\hat{d})\le \frac{32\De}{c^2}\cdot \diam\,(K,\hat{d})\,.$$\endp



We now introduce the last ingredient for the proof of Theorem \ref{prop:gluqm}.
Let $X_1$ and $X_2$ be two closed subsets of a metric
space $X$ such that $X_1\cap X_2\ne \emty$. The {\it seam} 
is by definition the closed set $Y=X_1\cap X_2$.
Following Agard and Gehring,
the {\it angle} $\angle (X_1,X_2)$ between $X_1$ and $X_2$ 
is by definition the supremum over all $c >0$ such that, for any $(x_1,x_2)\in X_1\times X_2$,
$$|x_1-x_2|\ge c\cdot\inf_{y\in Y}\{|x_1-y|+|y-x_2|\}\,.$$
We recall \cite[Thm\,3.1]{ph:sewing}:

\REFPROP{thflat} Let $X=X_1\cup X_2$ and $X'=X_1'\cup X_2'$ be 
metric spaces with positive angles. Let us assume that $Y=X_1\cap X_2$
and $Y'=X_1'\cap X_2'$ are $\la$-uniformly perfect subspaces such that
$\diam\, Y\ge\mu \diam X_1$ for some $\mu\in (0,1)$.

If $f:X\to X'$ is a \homeo such that $f|_{X_j}$ is $\eta$-quasisymmetric
and $f(X_j)=X_j'$, then $f$ is globally $\heta$-quasisymmetric quantitatively.\ENDPROP

We now need to check that angles remain positive when applying Lemma \ref{lma:stereo}:

\REFLEM{lma:positivestereo} For any $c> 0$ and $\De \ge 1$, there exists $\hat{c}>0$ with the following property.
Let $A$, $B$  be subsets of a metric space $(Z,d)$ with $A\cap B = X$ and
$\angle_d (A,B) \ge c $. Let $w\in Z$ and let us consider the metric space $(Z\setminus\{w\},\hat{d})$ given by
Lemma \ref{lma:stereo}. Then $\angle_{\hat{d}} (A\setminus\{w\}, B\setminus\{w\})\ge \hat{c}$.
\ENDLEM

\Proof  Fix $a\in A\setminus \{w\} $, $b\in B\setminus\{w\}$ and $x\in X$ such that $d(a,b)\ge c(d(a,x)+d(x,b))$.

We first observe that it is enough to find a constant $C=C(c,\De)$ and a point $x'\in X\setminus\{w\}$ such that
$$\min\{\hat{d}(a,x'),\hat{d}(b,x')\}\le C \hat{d}(a,b)\,.$$
In this case, it follows  from the triangle inequality that  $\max\{\hat{d}(a,x'),\hat{d}(b,x')\}\le (C+1)  \hat{d}(a,b)$ so that
$$\hat{d}(a,b)\ge \frac{1}{2(C+1)}(\hat{d}(a,x')+\hat{d}(x',b))\,.$$

We rely on the notations of Lemma \ref{lma:stereo} and shall use $q$ instead of $\hat{d}$.
Let $\al \in (0,1)$ be a small constant which will only depend on 
$\De$.
If $\de(x)\ge \al \de(a)$ then $$q(a,b)\ge c  \frac{d(b,x)}{\de(a)\de(b)}\ge  c\al  q(b,x)\,;$$
similarly if $\de(x)\ge \al \de(b)$.

Let us now assume that $\de(x)\le \al\min\{\de(a),\de(b)\}$. On the hand, the triangle inequality implies 
that $$d(a,b)\ge c d(b,x)\ge c(\de(b)-\de(x))\ge c(1-\al)\de(b)$$ holds
so that $q(a,b)\ge c(1-\al)/\de(a)$.
On the other hand, there is a point $x'\in X$ such that $$d(x,x')\ge \frac{1}{2\De} d(x,a)\ge \frac{1}{2\De} (1-\al)\de(a)$$ so that
$$\de(x')\ge d(x,x')-\de(x)\ge \left(\frac{1-\al}{2\De}-\al\right)\de(a)\,.$$
Choosing $\al$ small enough with respect to $\De$, we get $\de(x')\ge (1/4\De)\de(a)$, so that
$$q(a,x')\le \frac{\de(a)+\de(x')}{\de(a)\de(x')}\le \frac{1+4\De}{\de(a)}\le \frac{1+4\De}{c(1-\al)} q(a,b)\,.$$\endp

We may now turn to the proof of the main result of the section.

\Proof (Theorem \ref{prop:gluqm})  We first notice that Theorem \ref{prop:gsewing} and Proposition \ref{prop:qmvsqs} imply that there is a distortion
function $\te_1$ such that the maps
$(X,d_\Si)\stackrel{h_X}{\longrightarrow} (X,d_{\Si})$ and $(L_K,d_\Si)\stackrel{h_K}{\longrightarrow} (L_{h_X(K)},d_{\Si})$, $K\in\PPP$,  
are $\te_1$-quasi-M\"obius.

 Let us pick $w\in X\,(\subset\Si)$  and apply Lemma \ref{lma:stereo} to $(\Si,w)$ and denote by $(Z,d)$
the resulting metric space.  Let $(Z',d')$ be the metric space obtained from $(\Si,h(w))$. The theorem follows if we  prove
that $h:Z\to Z'$ is quasisymmetric with a distortion function which only depends on $\te$, $\eta$ and $\De$.

Proposition \ref{prop:qmvsqs} implies that $h_X: ( X\cap Z,d)\to ( X\cap Z',d')$ is $\eta_1$-quasisymmetric, where $\eta_1$
only depends on $\te$. Similarly, there is some $\te_2$ which only depends on $\te$ such that
$(L_K\cap Z,d)\stackrel{h_K}{\longrightarrow} (L_{h_X(K)}\cap Z',d')$ is $\te_2$-quasi-M\"obius for all $K\in\PPP$.

Now, from Theorem \ref{prop:gsewing}, 
we know the existence of $\De_1$ such that 
$\diam_\Si L_K\le \De_1\diam_\Si K$
for all $K\in\PPP$. Lemma \ref{lma:diamstereo} implies the same property on $Z$ and $Z'$ with another constant 
$\De_2$. So by Corollary \ref{cor:qmvsqs},
there is a distortion function $\eta_2$ such that  $(L_K\cap Z,d)\stackrel{h_K}{\longrightarrow} (L_{h_X(K)}\cap Z',d')$ is $\eta_2$-quasisymmetric 
for all $K\in\PPP$. 

For all $K\in\PPP$, we have $(L_K\cap X)= K$, and the separating
property (\ref{eq:0}) implies that  $\angle (L_K, X)\ge c_0$ for some constant which only depends
on $\eta$ and $\De_0$. Lemma \ref{lma:positivestereo} provides us with uniform positive angles in $Z$ and $Z'$. 
Therefore, Proposition \ref{thflat} implies that $h: (L_K\cup X)\cap Z\to (h(L_K)\cup X)\cap Z'$ is $\eta_3$-quasisymmetric
for some $\eta_3$ which only depends on $c_0$, $\eta_1$ and $\eta_2$, so on $\eta$ and $\De_0$. 
Similarly,  given $K_1\ne K_2\in\PPP$,  $(L_{K_1}\cup X)\cap (L_{K_2}\cup X)=X$ and   $\angle ((L_{K_1}\cup X), (L_{K_2}\cup X))\ge c_0$
so that $h|_{(L_{K_1}\cup L_{K_2}\cup X)\cap Z}$ is $\eta_4$-quasisymmetric onto its image, with
 $\eta_4$ which only depends on $\eta$ and $\De_0$. 
 If we pick a third $K_3\in\PPP$, then we obtain that $h|_{(L_{K_1}\cup L_{K_2}\cup L_{K_3}\cup X)\cap Z}$ is $\eta_5$-quasisymmetric onto its image, with
 $\eta_5$ which only depends on $\eta$ and $\De_0$. 

 This is enough to conclude that $h:Z\to Z'$ is $\eta_5$-quasisymmetric.\endp


\section{ Convergence group actions}\label{sec:hyp}
Let $G$ be a group acting by homeomorphisms on a metrizable compact  space.
The action of $G$ is {\em a  convergence action} if its diagonal action on the set of distinct triples is
properly discontinuous. We say the action is {\em uniform} if its action is also cocompact on
the set of distinct triples. As for Kleinian groups, the {\it limit set}  $\La_G$ is by definition
the unique minimal closed  invariant subset of $X$. It is empty if $G$ is finite.

We define word hyperbolic groups, their boundaries at infinity, their action on their boundary and quasiconvex subgroups of word hyperbolic
groups with the following characterization due to Bowditch \cite{bowditch:characterization,bowditch:convergence_groups}.
A more standard approach may be found for instance in \cite{gromov:hyperbolic,ghys:delaharpe:groupes,
kapovitch:benakli}.

\begin{theorem}[Bowditch]\label{thm:def:hyp} Let $G$ be a convergence group acting on a perfect metrizable space $X$.
The action of $G$ is uniform on $\La_G$ if and only if $G$ is word hyperbolic and there exists an equivariant homeomorphism
between $\La_G$ and the boundary at infinity $\partial G$ of $G$. Moreover, a subgroup $H$ is a quasiconvex subgroup of a word hyperbolic
$G$ if and only if $\La_H$ has at most two points or if the restriction of the action of $H$ to $\La_H$ is uniform.\end{theorem}


When $G$ is one-ended, then its
boundary  is connected, locally connected and without
(global) cut points \cite{bestvina:mess, swarup:cutpoint}. Note also that a word hyperbolic group is finitely presented, and is thus
{\it accessible} \cite{dunwoody:accessibility}: it can be decomposed over finite groups into an amalgamated free product of finitely many subgroups
so that each factor is a finite or one-ended.

Given a metric in its gauge $\GGG(G)$, the group $G$ acts as a uniform
quasi-M\"obius group \cite{paulin:determined}.
Moreover, there exists a constant $m>0$ such that, for any distinct points $x_1,x_2,x_3\in\partial G$,
there is some $g\in G$ such that $\{g(x_1),g(x_2),g(x_3)\}$ is $m$-separated. It follows that $\partial G$
is doubling and LLC when one-ended,   see  \cite{ph:bbki} and the references therein.

\subsection{Quasiconvex subgroups of word hyperbolic groups}
We establish some general properties in the spirit of  \cite[Thm\,5]{kapovich:kleiner:lowdim} and \cite[Prop.\,1.4]{bonk:sierpinski}
where the authors dealt with carpets. They are motivated by Proposition \ref{prop:cck}.

\REFPROP{prop:ginvns} Let $G$ be a non-elementary hyperbolic group and $\PPP$ a $G$-invariant
null-sequence in $\partial G$, each element containing at least two points. Then
\bit
\item[(a)] the set $\PPP/G$ is finite;
\item[(b)] there exists a distortion function $\eta$ such that, for any $K\in\PPP$ and $K'\in G(K)$,
there exists $g\in G$ such that $g|_K:K\to K'$ is $\eta$-quasisymmetric; 
\item[(c)] for any $K\in\PPP$, the stabilizer $G_K$ of $K$ is infinite
and acts on $K$ as a uniform convergence
group;
\item[(d)] the boundary $\partial G$ is doubling and porous relatively to $\PPP$ if the elements of $\PPP$ are connected;
\item[(e)] if the elements of $\PPP$ are pairwise disjoint and connected, then they are uniformly relatively separated i.e.,
 there is a constant $s>0$ such that $\dist(K_1,K_2)\ge s \min\{\diam\,K_1,\diam\,K_2\}$ for every distinct
pairs $K_1,K_2\in\PPP$.
\eit\ENDPROP

In other words, (c) means that $G_K$ is a {\it quasiconvex subgroup} of $G$.
These properties may be essentially established with the  conformal elevator principle
\cite[Prop.\,4.6]{ph:bbki}:

\REFPROP{prop:conf_elevator}{\em\bf (Conformal elevator principle)}
Let $G$ be a non-elementary hyperbolic group and consider its boundary $\partial G$ endowed with
a metric from its gauge. Then there exist  definite sizes $r_0\ge \de_0>0$
and a distortion function $\eta$ such that, for any $x\in X$, and any $r\in (0,\diam\,\partial G/2]$, 
there exists $g\in G$ such that
$g(B(x,r))\supset B(g(x),r_0)$, $\diam\,B(g(x),r_0)\ge 2\de_0$ and $g|_{B(x,r)}$ is $\eta$-quasisymmetric.
\ENDPROP

We draw the following consequence of the conformal elevator principle. Let $\ep\in (0,1)$, and let us consider
$y\in B(x,r)$ such that $d(g(x),g(y))\ge \de_0$, and let $z\notin B(x,\ep r)$. Then
either $z\notin B(x,r)$ so that $d(g(x),g(z))> r_0$, or 
$$d(g(x),g(z))\ge \frac{1}{\eta(1/\ep)} d(g(x),g(y))$$ so that we have
\begin{equation}\label{eq:subball} 
g(B(x,\ep r))\supset B(g(x), \de_0/\eta(1/\ep))\end{equation}

\Proof (Prop.\,\ref{prop:ginvns}) Let us fix a metric in $\GGG(G)$ and let $m>0$ be such that
any distinct triple of $\partial G$ can be $m$-separated by an element of $G$. 
Given $\de>0$, we let $\PPP_\de$ denote the subset 
of elements $K$ of $\PPP$ such that $\diam\,K\ge \de$; this set is  finite since $\PPP$ is a null-sequence
and non-empty for small enough $\de$. 

For all $K\in\PPP$, we can find two points 
$x_1,x_2\in K$ and a group element $g\in G$ such that $\{g(x_1),g(x_2)\}$ is $m$-separated: this implies
that $g(K)\in\PPP_m$, so that $\PPP$ is composed of finitely many orbits and (a) holds.


Let $r_0$ and $\de_0$ be the constants arising from the conformal elevator principle.
To prove (b), we notice that since $\PPP_{\de_0}$ is finite and that quasi-M\"obius mappings between compact
sets are quasisymmetric,
it is enough to prove that any $K\in\PPP$ can be mapped by a uniform quasisymmetric map to
an element of $\PPP_{\de_0}$. But this is exactly what the conformal elevator principle does.

Let us fix $K\in\PPP$ and let us first assume that $K$  contains at least three points. 
Let us enumerate $\PPP_m\cap G(K)=\{K_1,\ldots, K_N\}$.
For each $j\in\{1,\ldots, N\}$, one can find $g_j\in G$ such that $g_j(K_j)=K$. Since $\PPP_m\cap G(K)$ is finite
and $\{g_j^{-1},\, j=1,\ldots, N\}$ are  uniformly continuous, there is some $m'\in (0,m]$ such that, 
for all $j\in\{1,\ldots, N\}$
and all $x,y\in K_j$ with $d(x,y)\ge m$, one has $d(g_j(x),g_j(y))\ge m'$.  
Since $G_K$ is a subgroup of $G$, its action on the set of distinct triples of $K$ is automatically properly discontinuous.
Let us prove that it is also cocompact. Let $x_1,x_2,x_3\in K$ and consider
$g\in G$ such that $\{g(x_1),g(x_2),g(x_3)\}$ is $m$-separated. It follows that $g(K)=K_j$ for some $j$,
hence $(g_j\circ g)\in G_K$ and  $\{(g_j\circ g)(x_1),(g_j\circ g)(x_2),(g_j\circ g)(x_3)\}$ is $m'$-separated.

If $K$ has only two points $x,y$, then it suffices to prove that its stabilizer is infinite: this will imply it
is two-ended, hence is quasiconvex
in $G$. Pick a sequence $(x_n)$ accumulating $x$. As above, we may find infinitely many $g_n\in G$
such that $g_n\{x,x_n,y\}$ are $m$-separated; the same argument as above proves that $G_K$ is infinite. 
This proves (c).

Let us prove that $\partial G$ is relatively doubling. 
Pick $x\in \partial G$ and $r\in (0,\diam\,G/2]$ and let us apply Proposition \ref{prop:conf_elevator}:
 there exists $g\in G$ such that
$g(B(x,r))\supset B(g(x),r_0)$ and $g|_{B(x,r)}$ is $\eta$-quasisymmetric. 

Since the restriction of $g$ is $\eta$-quasisymmetric, it follows from Lemma \ref{lma:qsbounds} 
that, for each boundary component $K$ which intersects $B(x,r)$,
$$2\eta \left(\frac{\diam\,B(x,r)}{\diam\,(B(x,r)\cap K)}\right)\ge \frac{\diam\,g(B(x,r))}{\diam\,g(B(x,r)\cap K)}\,.$$
So if we assume that $\diam\,(B(x,r)\cap K) \ge \ep r$, then $\diam\,g(B(x,r)\cap K)\ge \de_0/\eta(2/\ep) $.
But since  $\PPP$ is a null-sequence, 
there are only finitely many such boundary components.

We now assume that the elements of $\PPP$ are connected.

The proof that $\partial G$ is relatively porous is similar: we apply the conformal elevator principle as above.
It follows from (\ref{eq:subball}) that $g(B(x,r/2))\supset B(g(x),r_0/\eta(2))$.
But since the action of $G$ is minimal on $\partial G$ which is compact, there exists a constant $\de_1>0$ such that,
for $y\in \partial G$, there exists $K_y\in\PPP$ with the property that $K_y\cap B(y,r_0/\eta(2))\ne\emty$ and 
$\diam\,K_y\ge \de_1$. Let $K_{g(x)}'$ be a connected component of $K_{g(x)}\cap g(B(x,r))$ which intersects $B(g(x),r_0/\eta(2))$;
it follows from  above that $\diam\,K_{g(x)}'\ge c r_0$ for some universal constant $c>0$. 
By construction, $K'= g^{-1}(K_{g(x)}')$ is a subset of an element of $\PPP$ which intersects $B(x,r/2)$;
moreover, Lemma \ref{lma:qsbounds} implies that $$\diam\,K' \ge (1/2\eta(2/c)) r\,.$$

We now turn to the proof of (e). Let $K_1$ and $K_2$ be two distinct boundary components such
that $\diam\,K_1\le \diam\,K_2$. 
Since $\PPP$ is a null-sequence, we may assume that $\diam\,K_1\le r_0$, and that $\dist(K_1,K_2)\le \diam\,K_1/2$.
Let $x\in K_1$ and $y\in K_2$ be such that $d(x,y)=\dist (K_1,K_2)$ and set $r=\diam\,K_1$. Let $K_2'$ be the component  
of $K_2\cap B(x,r)$ which contains $y$ so that $\diam\,K_2'\ge r/2$ and $\dist(K_1,K_2')=\dist(K_1,K_2)=d(x,y)$.
Apply the conformal elevator
principle to $B(x,r)$. It follows that $\diam\,g(K_1)\ge r_0$ and, from Lemma \ref{lma:qsbounds}, we may deduce
that $$\diam\,g(K_2)\ge \diam\,g(K_2')\ge \frac{r_0}{2\eta(2)}\,.$$
Since $\PPP$ is a null-sequence and the components are pairwise disjoint, there is a constant $s_0>0$ independent
from $K_1$ and $K_2$ such that $$d(g(x),g(y))\ge \dist(g(K_1), g(K_2'))\ge  \dist(g(K_1), g(K_2)) \ge s_0 r_0\,.$$
Hence, applying a last time Lemma \ref{lma:qsbounds} yields 
$$d(x,y)\ge \frac{ r}{2\eta^{-1}(s_0)}\ge\frac{\diam\,K_1}{2\eta^{-1}(s_0)} $$ and so 
$$\dist(K_1,K_2)\ge \frac{1}{2\eta^{-1}(s_0)}\min\{\diam\,K_1,\diam\,K_2\}\,.$$ \endp

\subsection{The extension property}
We give a criterion which enables us to extend a convergence action to a larger space.
We will base our approach on the following notion.

\begin{defn}[Extension property] Let $X$ be a metrizable continuum and $Y\subset X$ be a compact subset.
We say that the pair $(Y,X)$ has the extension property 
if any convergence action of a group on $Y$ is the restriction of a convergence action on $X$ by
the same group. If $X$ is supplied with a metric, we  say that the pair $(Y,X)$ has the conformal 
extension property 
if any action of a group which acts
by uniform quasi-M\"obius homeomorphisms on $Y$ can be extended to an action on $X$ by uniform quasi-M\"obius 
homeomorphisms.\end{defn}

\begin{rmk}{\rm A discrete group $G$ of uniformly quasi-M\"obius self-homeomorphisms 
of a compact metric space $Z$  has the convergence property.
}
\end{rmk}

These extension properties are motivated by the following result due to Casson and Jungreis \cite{casson:jungreis},
Gabai \cite{gabai:S1}, 
 Hinkkanen 
\cite{hinkkanen:uniform_qsg, hinkkanen:structure_qsg}
and Markovic \cite{markovic:quasisymmetric:groups}.

\REFTHM{thm:qmonS1}
\indent (a) Any faithful convergence action of a group on the unit circle is conjugate to an action of
a Fuchsian group.

(b) Any uniformly quasi-M\"obius group of homeomorphisms on the unit circle
 is quasisymmetrically conjugate to a group of M\"obius transformations.\ENDTHM
 
So we may conclude that the extension property is not void. 
\REFCOR{cor:qmonS1} The pair $(\SS^1,\DD)$ has both the extension and conformal extension
properties.\ENDCOR

\Proof The extension property is a direct consequence of Theorem \ref{thm:qmonS1} (a):
if $G$ is a convergence group of $\SS^1$, there is a 
homeomorphism $h:\SS^1\to \SS^1$ such that $G'=hGh^{-1}$ is a group of M\"obius
transformations. Therefore, this group $G'$ acts canonically on $\DD$ by
M\"obius transformations as well. 
Let $H:\DD\to \DD$ be a homeomorphism  which extends $h$.
The group $H^{-1}G' H$  extends the action
of $G$ faithfully. 

Assuming $G$ is a group of uniform quasi-M\"obius homeomorphisms, Theorem \ref{thm:qmonS1} (b) provides us with  a quasisymmetric homeomorphism $h:\SS^1\to \SS^1$ such that $G'=hGh^{-1}$ is a group of M\"obius
transformations. 
Let $H:\DD\to \DD$ be a quasiconformal map which extends $h$ \cite{ahlfors:lectures_qc}.
The group $H^{-1}G' H$ is a uniform group of quasi-M\"obius maps which extends the action
of $G$ faithfully.\endp

\bigskip

\REFTHM{thm:extcva} Let $G$ be a group acting on a metrizable continuum $X$ as a convergence group action.
Let $K_0\subset X$ be a subcontinuum such that $G(K_0)$ defines an admissible sequence.
We also assume that there exists an embedding $\psi:K_0\to L$, where $L$ is a metrizable  
continuum such that $(\psi(K_0),L)$ has the extension property. 
For each $K\in G(K_0)$, let $g_K\in G$ map $K_0$ to $K$, and let us consider the space
$\Si$ obtained by Proposition \ref{prop:sicomp} with gluing maps $(\psi\circ g_K^{-1})$.
Then we may extend the action of $G$ on $X$ to $\Si$ as a convergence group action  with same
limit set.\ENDTHM

\Proof Set $K'=\psi(K_0)$, let $\vp_0:K'\to K_0$ denote its inverse and let $\vp_K=g_K\circ \vp_0:K'\to K$ 
and $\psi_K=\vp_K^{-1}$.
For each $K$, we set $L_K= L$ and we let $\Si$  be the 
metric space obtained from Proposition \ref{prop:sicomp}.

Denote by $H$ the set of self-homeomorphisms of $K'$ of the form $\vp_{g(K)}^{-1}\circ g\circ \vp_K$, among
all $K\in G(K_0)$ and $g\in G$. Given $h_1$ and $h_2$ in $H$, we let $g_1,g_2\in G$ and $K_1,K_2\in G(K_0)$ be
such that $h_j= \vp_{g_j(K_j)}^{-1}\circ g_j\circ \vp_{K_j}$.
Then 
\begin{eqnarray*}
h_1\circ h_2 & = & \vp_{g_1(K_1)}^{-1}\circ g_1\circ \vp_{K_1}\circ \vp_{g_2(K_2)}^{-1}\circ g_2\circ \vp_{K_2}\\
& = & \vp_{g_1(K_1)}^{-1}\circ g_1\circ (g_{K_1}\circ \vp_0)\circ (\vp_0^{-1}\circ g_{g_2(K_2)}^{-1})\circ g_2\circ \vp_{K_2}\\
& = & \vp_{g_1(K_1)}^{-1}\circ  ( g_1\circ g_{K_1}\circ g_{g_2(K_2)}^{-1}\circ g_2)\circ \vp_{K_2}\,.
\end{eqnarray*}
We may check that $( g_1\circ g_{K_1}\circ g_{g_2(K_2)}^{-1}\circ g_2) (K_2)= g_1(K_1)$ so that $h_1\circ h_2\in H$
and $H$ is a group of homeomorphisms. Note that $H$ is isomorphic to the stabilizer of $K_0$:
if $h=\vp_{g(K)}^{-1}\circ g\circ \vp_K$, then we may write $$h = \vp_{0}^{-1} \circ (g_{g(K)}^{-1}\circ g\circ g_K)\circ\vp_0$$
and $$ (g_{g(K)}^{-1}\circ g\circ g_K)(K_0)=(g_{g(K)}^{-1}\circ g)(K)=g_{g(K)}^{-1}(g(K))=K_0\,.$$
Therefore, $H$ is  a convergence group on $K'$.

By the extension property, $H$
acts on $L$ as a convergence group as well.

We now define an action of $G$ on $\Si$. 
Let $g\in G$.  
Fix $K\in G(K_0)$, and let $h=\vp_{g(K)}^{-1}\circ g\circ \vp_K:K'\to K'$ and $\hat{h}:L_K\to L_{g(K)}$
be its extension. 
By construction,
$g|_K= (\psi_{g(K)}^{-1}\circ \hat{h}\circ\psi_K)|_K$ so that Lemma \ref{lma:homeoglobal} implies that
these maps patch up into a homeomorphism $\hat{g}:\Si\to\Si$.

Let $\hat{g}_j\in G$, $j=1,2$. On $X$, we find $g_1,g_2$ such that
$\hat{g}_j=g_j$ so that $(\hat{g}_1\circ \hat{g_2})|_{X}= (g_1\circ g_2) |_{X}$. Set $g=g_1\circ g_2$
and $\hat{g}$ its extension to $\Si$. We have to prove that $\hat{g}_1\circ \hat{g_2}=\hat{g}$.

Fix $K\in G(K_0)$, and let us prove that $(\hat{g}_1\circ \hat{g_2})|_{L_K}= \hat{g} |_{L_K}$.

We let $(h_1,\hat{h}_1)$, $(h_2,\hat{h}_2)$ and $(h,\hat{h})$  be associated to $(g_2,K)$,  $(g_1, g_2(K))$, and $(g,K)$.  
On $K'$, we have
\begin{eqnarray*}
h_1\circ h_2 & = &  ( \vp_{g_1(g_2(K))}^{-1}\circ g_1\circ \vp_{g_2(K)})\circ (\vp_{g_2(K)}^{-1}\circ g_2\circ \vp_K)\\
& = & \vp_{(g_1\circ g_2)(K)}^{-1}\circ (g_1\circ  g_2)\circ \vp_K\\
& = & \vp_{g(K)}^{-1}\circ g\circ \vp_K= h\,.\end{eqnarray*}
Hence, the extension property implies that $\hat{h}_1\circ \hat{h}_2=\hat{h}$.

It follows that the extended maps define an action of $G$ on $\Si$. Since $(G(K_0))$
is admissible, the limit set $\La$ of this action coincides with its embedded copy in  $X$. This embedding
is  equivariant by construction. 
\endp

We now apply this construction to hyperbolic groups.

Let $G$ be a one-ended hyperbolic group, and let us supply $\partial G$ 
with a metric from its gauge. We assume that there exists a continuum $K_0\subset \partial G$
such that $G(K_0)$ forms an admissible sequence of subcontinua of $\partial G$. We also
assume that there exists a quasisymmetric embedding $\psi:K_0\to L$, where $L$ is a metric 
continuum such that $(\psi(K_0),L)$ has the conformal extension property. 

\REFTHM{thm:extpr} There exist gluing functions of $L$  along $G(K_0)$ 
so that the compact metric space $\Si$ constructed by Theorem \ref{prop:gsewing} 
can be endowed with an action of $G$ by uniformly quasi-M\"obius maps such that
there is an equivariant bi-Lipschitz homeomorphism $\vp: \partial G \to \La$
where $\La\subset\Si$ is the limit set of the action of $G$.
\ENDTHM

\Proof
According to Proposition \ref{prop:ginvns} (b), there exists $\eta$
such that, for any $K\in G(K_0)$, there is some $g_K:K_0\to K$ which is $\eta$-quasisymmetric.
We use the same notation as above and note that $\diam\,L_K\le \De_0 \diam\, \psi_K(K)$ with $\De_0=\diam\,L/\diam\,K'$,
and that $\vp_K$ and $\psi_K$ are uniformly quasisymmetric.
We let $\Si$  be the 
metric space obtained from Proposition \ref{prop:sicomp} and Theorem \ref{prop:gsewing}.
 The embedding of $\partial G$ in $\Si$ is bi-Lipschitz by Theorem \ref{prop:gsewing}.

Since $G$ is uniformly quasi-M\"obius and the $\vp_K$'s are uniformly
quasisymmetric, the group $H$ defined as in Theorem \ref{thm:extcva} is also a uniform quasi-M\"obius group.
According to Theorem \ref{thm:extcva}, the group $G$ acts as a convergence group on $\Si$ with limit set $\partial G$.
Theorem \ref{prop:gluqm} implies that this action is also uniformly quasi-M\"obius.

This embedding
is bi-Lipschitz by Theorem \ref{prop:gsewing} and equivariant by construction. 
\endp

\subsection{Planar actions of word hyperbolic groups}\label{subsec:planaraction}
For the proof of Theorem \ref{thm:main}, we need some properties on 
 the boundaries of hyperbolic groups. Some may be established just using the 
planarity. Others seem more specific to planar actions.
Let us recall that a space is {\it planar} if it can be embedded topologically in the two-sphere $S^2$, and
a word hyperbolic group admits a {\it planar action} if its boundary admits an embedding in $S^2$
in such a way that the action of each element of the group is the restriction of a self-homeomorphism of $S^2$.

We start with a lemma which describes the topology of the boundary of a one-ended hyperbolic group provided it is planar.

\REFPROP{prop:planarboundary} Let $G$ be a one-ended hyperbolic group with a planar boundary.
Then either $\partial G$ is homeomorphic to the sphere or $\partial G$ is a degenerate carpet. 
In the latter case, any embedding defines  an admissible collection
of  boundary components.\ENDPROP

\Proof If $\partial G$ is not a sphere, then it is one-dimensional. But the boundary
of a one-ended group has no global cut point and is locally connected  \cite{bestvina:mess, swarup:cutpoint}. Therefore, Lemma \ref{lma:degcarpet}
implies that it is a degenerate carpet.

Let us fix an embedding $\vp:\partial G\to S^2$ and define $\PPP$ to be the collection of
boundaries of the components of $S^2\setminus\vp(\partial G)$. 
Since $\partial G$ is locally connected, it follows from  \cite[Thm.\,VI.4.4]{whyburn:analytic_topology} 
that $\PPP$ is an admissible null-sequence. \endp


We may now establish some stronger assumptions when $G$ admits a planar action.

\REFPROP{prop:sierpgroup} Let $G$ be a one-ended hyperbolic group with a planar action.
Then
\ben
\item the boundary components form an invariant null-sequence;
\item every boundary component $K$ is a uniform quasicircle and $G_K$ is virtually a cocompact Fuchsian group.
\een
\ENDPROP

\Proof (Prop.\,\ref{prop:sierpgroup}) Let $\vp:\partial G\to S^2$ be the topological embedding such that
$G$ leaves invariant the set $\PPP$ of boundary components of $S^2\setminus\vp(\partial G)$.  
We know from Proposition \ref{prop:planarboundary} that $\PPP$ is  also a null-sequence.
Proposition \ref{prop:planarboundary} also implies that each boundary component
is a Jordan curve. 

By Proposition \ref{prop:ginvns}, for any $K\in\PPP$, $G_K$ is a uniform convergence group 
on $K$. 
We may then conclude from  Theorem \ref{thm:qmonS1} that $G_K$ is virtually isomorphic to
a cocompact Fuchsian group. Therefore, $K$, as the limit set of a quasiconvex subgroup, is doubling and LLC  as explained in the introduction of \S\,\ref{sec:hyp}; hence
$K$ is a quasicircle  by \cite{tukia:vaisala:qs}. Since $G$ acts by uniform quasi-M\"obius maps and there are
only finitely many orbits of boundary components, we conclude that $\PPP$ is formed of uniform quasicircles.
\endp


 



We may now prove Theorem \ref{thm:main}.

\Proof (\Th \ref{thm:main})  The necessity follows from \cite{sullivan:ihes50} as mentioned in the introduction. 

Let $\vp:\partial G\to S^2$ be the topological embedding such that
$G$ leaves invariant the set $\PPP$ of boundary components of $S^2\setminus\vp(\partial G)$. 
By Proposition \ref{prop:sierpgroup}, each peripheral circle
is uniformly quasisymmetric equivalent to the unit circle.

Note that $\partial G$ is LLC, relatively doubling and
porous with respect to its boundary components. 
Therefore, Corollary \ref{cor:qsembedcarpets}
applies, so $\partial G$ can be
quasisymmetrically embedded in the sphere. 

Since $(\SS^1,\DD)$ has the extension property according to Corollary \ref{cor:qmonS1},
Theorem \ref{thm:extpr} enables us to extend the action of $G$ to the whole sphere as
a group of uniform quasi-M\"obius transformations. 
By Sullivan's straightening theorem (Theorem \ref{thm:sullivan}), 
this action is conjugate to a group of M\"obius transformations.
The action being properly discontinuous on the complement of $\partial G$, the group is
discrete. 

If $G$ does not act faithfully on $\partial G$, then there is a normal finite subgroup $F$
such that $G/F$ is isomorphic to a discrete subgroup of isometries of $\HH^3$.
In order to conclude, we use material developed in  \S\S\,\ref{sec:special} and \ref{sec:qcerf}.
According to Corollary  \ref{prop:virtspeck} below and  \cite[Lemma 2.10]{agol:virtualHaken}, the group $G$ 
is virtually special, hence virtually torsion-free by Corollary \ref{cor:virtorsionfree}.

Therefore, $G$ has a torsion-free finite index subgroup which is isomorphic to a convex-cocompact Kleinian group.
\endp

\begin{rmk} {\rm In the case of carpets, one could construct the Kleinian group differently.
From \cite[Thm 1.1]{bonk:sierpinski},  we know that  $\partial G$ is quasisymmetric
equivalent to a round carpet $\La$ (with all boundary components a round circle)
of measure $0$ (since the carpet is porous in $\cbar$), and $G$ acts by quasisymmetric homeomorphisms. 
Therefore, \cite[Thm 1.1]{bonk:kleiner:merenkov} implies that
each element of $G$ acting on $\La$ is actually the restriction of a M\"obius transformation.}
\end{rmk}

\Proof (Remark \ref{rmk:gdcboundary})
There are several approaches to see whether $G$ is quasi-isometric to a convex subset of $\HH^3$
with geodesic boundary or not. The first is based on hyperbolic geometry. If $G$ is quasi-isometric
to a convex subset with geodesic boundary, then it is virtually the fundamental group of a compact hyperbolic
manifold with geodesic boundary. It follows that this manifold is acylindrical so the limit set is
homeomorphic to the Sierpi\'nski carpet, see e.g. \cite{ctm:classification}. The converse is one of the
main steps of Thurston's hyperbolization theorem, see \cite{ctm:iteration}.

Another  route goes as follows.
When $\partial G$ is a carpet, peripheral circles are disjoint uniform quasicircles so they 
are also uniformly separated by Proposition \ref{prop:ginvns}. Hence, \cite[Thm 1.1]{bonk:sierpinski} 
implies that $\partial G$ is quasisymmetric
equivalent to a round carpet $\La$ (with all boundary component a round circle). 
It now follows from \cite{bonk:schramm:embeddings} that $G$ is quasi-isometric
to a convex subset of $\HH^3$ with geodesic boundary.



When $\partial G$ is not a carpet, then $\partial G$ admits local cut points \cite[Thm\,4]{kapovich:kleiner:lowdim}.
Either it is homeomorphic to a circle so that $G$ is virtually a cocompact Fuchsian group. Or the local
cut points are structured in equivalence classes as described by Bowditch \cite{bowditch:jsj}, cf.\,\S\,\ref{sec:jsj} below.
Note that local cut points have to belong to boundary components: 
if $x\in\partial G$ does not belong to any boundary component,
\cite[VI.4.5]{whyburn:analytic_topology} asserts that we can find a nested sequence of Jordan curves contained
in $\partial G$ enclosing $x$ in any neighborhood of $x$. This prevents $x$ to be a local cut point.
These cut points are associated to two-ended hyperbolic groups: every component of the complement which
contains one of the fixed points has to contain the other one: therefore, there are components which
intersect in at least two different points, preventing them to be mapped simultaneaously to round disks.\endp


We may now prove Corollary \ref{cor:planaraction}.

\Proof (Corollary \ref{cor:planaraction}) 
If $\La_G$ is connected, then  Theorem \ref{thm:main} and its proof imply that $G$ is isomorphic to a convex-cocompact
Kleinian group with conjugate actions by  a homeomorphism which extends to the whole sphere. 
Moreover, according to \cite[Cor.\,4.5]{martin:tukia:msri}, the action of each stabilizer of each component
is essentially unique, so we may build a conjugacy using similar  considerations as in the proof of Theorem
\ref{thm:extcva}.

Let us assume that $\La_G$ is not connected, and let us write $\oo_G=S^2\setminus \La_G$.
 We apply the decomposition techniques from \cite{abikoff:maskit}, see also \cite{martin:skora}.
Let $p:\oo_G\to\oo_G/G$ be the covering and $N$ be the normal subgroup defining the covering;
we consider a 
finite set of simple, disjoint loops $M=\{u_1,\ldots, u_n\}$ on $\oo_G/G$ such that 
 there exist (minimal) positive integers $a_1,\ldots, a_n$ so that the normal subgroup $N_M$ generated by
 $u_1^{a_1}, \ldots, u_n^{a_n}$ is a subgroup of  $N$. Let $\G=p^{-1}(\{u_1,\ldots, u_n\})$: this is a countable
set of pairwise disjoint homotopically non-trivial simple loops in $\oo_G$. 

We construct a tree $T$ as follows. Let the set of vertices be the connected components of
$S^2\setminus \cup_{\g\in\G}\g$ and put an edge between two such components if they share a
curve of $\G$ on their boundary. Let us observe that $\G$ is a null-sequence since $\G/G$ is finite
and the curves live in the set of discontinuity of $G$. 
Since each curve is a Jordan curve which separates $S^2$, it follows that
$T$ is simply connected and that the ends of $T$ correspond to a nested sequence of disks so
that $T$ is connected. This implies that $T$ is a tree.

By construction, $G$ acts simplicially on $T$ and $T/G$ is finite since $\G/G$ is finite. 
Moreover, an edge stabilizer corresponds to the stabilizer of a curve $\g\in\G$; since $\g\subset\oo_G$ and
$G$ is torsion-free, each edge stabilizer is trivial and there are no edge inversions. 
It follows that this action yields a decomposition
of $G$ as a free product. Since $G$ is hyperbolic, it is accessible so the number $n$ of edges in $T/G$ is bounded. 
But, according to \cite[Lemma 5]{maskit:planarity}, if $N_M\ne N$, then we can find $u_{n+1}$ disjoint from $M$ and $a_{n+1}$ 
with $u_{n+1}^{a_{n+1}}\in N$, so we may repeat the above construction with $M\cup\{u_{n+1}\}$. The accessibility of $G$ implies that this process has to stop, meaning that
we end up with a multicurve  $M$ such that $N_M=N$. 

This implies in particular that $\oo_G/G$ is a finite union of compact surfaces.
Let $T_0'\subset T/\G$ be a maximal tree and let us consider a connected lift  
$T_0\subset T$ so that each vertex is represented in $T$ exactly once; the stabilizers of vertices of $T_0$
are quasiconvex subgroups according to \cite[Prop.\,1.2]{bowditch:jsj}. Moreover, they are either trivial
or  one-ended for otherwise $\G/G$ would not generate $N$. 

Let $v\in T_0$; if its stabilizer $G_v$ is trivial, we may associate a $3$-ball
$M_v$; otherwise $G_v$ is a one-ended torsion-free planar group hence it is conjugate to a discrete group
of M\"obius transformations by Theorem \ref{thm:main}, and so $G_v$ is isomorphic to the fundamental
group of a hyperbolic manifold $M_v$. Each edge orbit attached to $v$ corresponds to a simple closed curve
on $\partial M_v$ which bounds a disk. The graph $T/G$ tells us how to build a manifold
$M$ by gluing the different $M_v$'s along the disks bounded by the curves above \cite{scott:wall}. 
We obtain in this way a Haken manifold which satisfies the assumptions of Theorem \ref{thm:thurston}.
Therefore, $G$ is isomorphic to a discrete subgroup $H$ of isometries of $\HH^3$. 

By the construction of $M$, the isomorphism between $G$ and $H$ yields an equivariant homeomorphism between
$\oo_G$ and $\oo_H$, so we may find a function $f:S^2\to \cbar$  such that (a) $f\circ G= H\circ f$,
(b) $f|_{\La_G}:\La_G\to\La_H$  and   $f|_{\Omega_G}:\oo_G\to\oo_H$ are  homeomorphisms.
Therefore, \cite[Prop.\,5.5]{bowditch:convergence_groups} implies that this conjugacy is a global
homeomorphism.\endp

\begin{rmk}\label{rmk:ahlforsfiniteness} The same argument as above proves that a convergence group $G$  on $S^2$ which is uniform on its limit
set satisfies the Ahlfors finiteness property: the orbit space $\oo_G/G$ is a finite union of closed surfaces (with finitely many
ramifications), cf. \cite{martin:skora}.\end{rmk}


\section{The JSJ decomposition}\label{sec:jsj}
We first summarize  briefly the JSJ-decomposition of a non-Fuchsian one--ended hyperbolic group $G$ following 
Bowditch \cite{bowditch:jsj}. 
Then we will focus on specific decompositions for groups with planar boundaries.

\subsection{General properties}
There exists a canonical simplicial minimal action of $G$ on a simplicial tree $T=(V,E)$ without edge inversions such that 
$T/G$ is a finite graph and which enjoys the following properties, cf. \cite[Thm\,5.28, Prop.\,5.29]{bowditch:jsj}.
If $v$ is a vertex (resp. $e$ an edge), we will denote by $G_v$ (resp. $G_e$) its stablizer, 
and by $\La_v$ (resp. $\La_e$) the limit set of $G_v$ (resp. $G_e$). Let $E_v$ denote the set of edges incident to $v\in T$.
Every vertex and edge group is  quasiconvex in $G$. 
Each edge group $G_e$ is two-ended and $\partial G\setminus \La_e$ is not connected. 
A  vertex $v$  of $T$ belongs to exactly one of the following three exclusive types.
\bit
\item[] {\bf Type I (elementary).---}  The vertex has bounded valence in $T$.  
Its stabilizer $G_v$  is  two-ended, and the connected components of  $\partial G\setminus \La_v$ are
in bijection with the edges incident to $v$. 
\item[] {\bf Type II (surface).---} The limit set $\La_v$ is cyclically separating and 
the stabilizer $G_v$ of such a vertex $v$ 
is a non-elementary virtually free group canonically isomorphic to a convex-cocompact Fuchsian group.
The incident edges are in bijection with the peripheral subgroups of that Fuchsian group. 
\item[] {\bf Type III (rigid).---} Such a vertex $v$ does not belong to a class above.  
Its stabilizer $G_v$ is a non-elementary quasiconvex subgroup. Every local cut point 
of $\partial G$ in $\La_v$ is in the limit set of an edge stabilizer incident to $v$; see Lemma \ref{lma:proptypeIII} 
for more properties of rigid type vertices.
\eit
No two vertices of the same type are adjacent, nor surface type and rigid can be adjacent either. 
The action of $G$ preserves the types. Therefore, the edges incident to a vertex $v$ 
of surface type or rigid type are split into finitely many $G_v$-orbits. 

Let $v$ be a vertex of $T$ of rigid type; if $e\in E_v$ is incident to $v$, let $C_v(e)$ denote the connected component of 
$\partial G\setminus \La_e$ which contains $\La_v$ and set $Z_e=\partial G\setminus C_v(e)$: 
this is a connected, locally connected, 
compact  set by construction. Let us define the following equivalence relation $\sim_v$ on $\partial G$. 
Say $x\sim_v y$ if $x=y$ or if there exists an edge $e\in E_v$ incident
to $v$ such that $\{x,y\}$ is a subset of  $Z_e$. 
Let $Q_v=\partial G/\sim_v$ be endowed with the quotient topology and let $p_v:\partial G\to Q_v$ be the canonical projection. 
Note that preimages of points are either
points or one of the $Z_e$'s, so they are always
connected and the map $p_v$ is {\it monotone}.

\REFLEM{lma:proptypeIII} If $v$ is of rigid type, then the space $Q_v$ is a  Hausdorff, compact, connected and locally connected space.
Moreover, no pair of points can disconnect $Q_v$ and
the local cut points of $Q_v$ correspond to the non-trivial classes of $\sim_v$ which 
disconnect locally $Q_v$ into exactly two components.
The group $G_v$ acts on $Q_v$ as a geometrically finite group and there are finitely many orbits of parabolic points.
\ENDLEM

Following Bowditch \cite{bowditch:rhg}, we say that a subgroup $H<G_v$ 
is {\it parabolic} if it is infinite,
fixes some point of $Q_v$, and contains no loxodromics. In this
case, the fixed point of $H$ is unique. We refer to it as a {\it parabolic point}.
The stabilizer  of a parabolic point is necessarily a parabolic group. There
is thus a natural bijective correspondence between parabolic points in
$Q_v$ and maximal parabolic subgroups of $G_v$. We say that a parabolic
group, $H$, with fixed point $p$, is {\it bounded} if the quotient $(Q_v \setminus \{p\})/H$
is compact. (It is necessarily Hausdorff.) We say that $p$ is a {\it bounded
parabolic point } if its stabilizer is bounded. A {\it conical limit point} is a point
$y\in Q_v$ such that there exists a sequence $(g_j)_{j\ge 0}$ in $G_v$, and distinct
points $a, b\in Q_v$, such that $g_j(y)$ tends to $a$ and $g_j(x)$ tends to $b$ for all $x \in Q_v\setminus \{y\}$.
We finally say that the action of $G_v$ on $Q_v$ is {\it geometrically finite} if every point is either conical or bounded parabolic
(they cannot be both simultaneously).

The proof of this  lemma will  rely on technical lemmas established by Bowditch in \cite{bowditch:jsj}.
 
\Proof 
Since $\partial G$ is locally connected, it follows that $\{Z_e\}_{e\in E_v}$ is a null-sequence, 
so  $\sim_v$ defines an upper semicontinuous decomposition of $\partial G$ and
$Q_v$ is Hausdorff (see also \cite[Cor.\,6.16]{carrasco:thesis} for a detailed proof based on the splitting). 
It follows that $Q_v$ is connected, locally connected and compact as the image 
under a continuous map of a Hausdorff,
connected, locally connected and compact set into a Hausdorff space. See \cite[Chapter VII, \S\,2,3]{whyburn:analytic_topology}
for details.

Cut points and local cut points of $Q_v$ yield (local) cut points of $\partial G$ by pull-backs under $p_v$. 
Let  $x\in Q_v$; if $x=p_v(Z_e)$, then it does not disconnect $Q_v$
since $Z_e$ does not ($\partial G\setminus Z_e=C_v(e)$); otherwise, since $v$ is of rigid type, 
it follows that $p_v^{-1}(x)$ is not a local cut point, so neither is $x$. 
We also conclude that the only possible  local cut  points correspond to the non-trivial classes. 

Fix an incident edge $e$. By construction, $\partial G\setminus Z_e$ has exactly two ends, each accumulating
one single point of $\La_e$. Hence, we may consider two disjoint connected
neiborhoods $N_1$ and $N_2$ of $\La_e$ in $(\partial G\setminus Z_e)\cup \La_e$; 
if there are small enough then $p_v(N_1)$ and $p_v(N_2)$ cover a neigborhood 
of $p_v(\La_e)$ and they intersect exactly at that point. This implies that it is a local cut point with two ends. 

The action of $G_v$ permutes the fibers of $p_v$, hence $G_v$ acts on $Q_v$. 
Since the action of $G_v$ is a convergence action and is minimal on $\La_v$ it is also the case on $Q_v$. 

If $e\in E_v$ is an incident edge, then $(\La_v\setminus \La_e)/G_e$ is compact by \cite{bowditch:convergence_groups},
so $(Q_v\setminus p_v(\La_e))/G_e= p_v(\La_v\setminus p_v(\La_e))/G_e$ is compact as well.
Thus, $p_v(\La_e)$ is a bounded parabolic point. 

Note that the action on $\La_v$ is uniform so every point is conical \cite{bowditch:convergence_groups}.  
Thus, if $y\in\La_v$, we may find a sequence $(g_j)_{j\ge 0}$ in $G_v$, and distinct
points $a, b\in \La_v$, 
such that $g_j(y)$ tends to $a$ and $g_j(x)$ tends to $b$ for all $x \in\La_v\setminus \{y\}$.
To conclude that $p_v(y)$ is also conical, it suffices to make sure that 
we may choose $a$ and $b$ not simultaneously in some $\La_e$.
Let us assume that $y$ is in no limit set of an incident edge and that  
there is indeed some incident edge $e$ such that $\La_e=\{a,b\}$. Pick a compact fundamental domain $K$
of $(\La_v\setminus \La_e)/G_e$. Then for any $j$, we may find $h_j\in G_e$ such that $h_jg_j(y)\in K$. 
Then it follows that, extracting a subsequence if necessary, $(h_j g_j(x))_j$ tends to $b$ 
for all $x\in\La_v$ but $(h_jg_j(y))_n$ remains far from $\La_e$. Therefore $p_v(y)$ is conical as well.

This proves that the action of $G_v$ on $Q_v$ is geometrically finite. 
Since there are only finitely many $G_v$-orbit of incident edges
to $v$ in $T$, there are only finitely many parabolic orbits in $Q_v$.

Let us prove that no pair of points separates $Q_v$ by contradiction. 
Following Bowditch \cite[\S\,3]{bowditch:jsj}, say two points $x$ and $y$ in $Q_v$ are equivalent, $x\sim y$,
if $x=y$ or $Q_v\setminus\{x,y\}$ is disconnected. Since every local cut point disconnects $Q_v$ locally into 
two components, this defines an equivalence relation on $Q_v$ according
to \cite[Lma\,3.1]{bowditch:jsj}. Moreover, \cite[Lma\,3.7]{bowditch:jsj} implies that each class is closed.
So, let us assume that $x\in Q_v$ belongs to a non-trivial class and let $y\sim x$, $y\ne x$. Since $x$ is necessarily
a local cut point, it is also parabolic. Hence, one can find a sequence $(g_n)$ stabilizing $x$ such that
$g_n(y)$ tends to $x$. It follows that no point is isolated in a non-trivial class. Therefore, each non-trivial
class is a perfect compact subset of local cut points of $Q_v$. But such a set is always uncountable, contradicting
that there are at most countably many parabolic points. Hence each class is trivial and no pair of points can
separate $Q_v$.
\endp

\subsection{Planar action of stabilizers of vertices of rigid type}
In this paragraph, we assume that $G$ is a one-ended hyperbolic group with planar boundary not homeomorphic to a simple closed curve. 
The goal is to analyze vertices of rigid type and interpret the incident edges as an {\it acylindrical paring} for the Kleinian manifold
when this vertex stabilizer is isomorphic to the fundamental group of a compact $3$-manifold with boundary.

Let us explain this terminology; more details can be found for instance in  \cite{thurston:bams, morgan:thurston, kapovich:book}. Let $M$ be a compact $3$-manifold with boundary. A surface $S$ is {\it properly embedded}
in $M$ if $S$ is compact and orientable and if either $S\cap\partial M=\partial S$ or
$S$ is contained in $\partial M$. A properly embedded surface $S$ is {\it incompressible} if $S$
is not homeomorphic to the $2$-sphere and if the inclusion $i:S\to M$ gives rise to an injective
morphism $i_*:\pi_1(S,x)\to\pi_1(M,x)$, for some base point $x\in S$. It turns out that a manifold is  Haken if it  contains
an incompressible surface. 
We say that $M$ has {\it incompressible boundary} if each component of $\partial M$ is incompressible. When $M$ is hyperbolizable, this is 
equivalent to the connectedness of the limit set of the group uniformizing $M$.
A compact {\it pared manifold} $(M,P)$ is given by a $3$-manifold $M$ as above together with a finite
collection of pairwise disjoint incompressible annuli $P\subset\partial M$ such that any cylinder in $C\subset M$
with boundaries in $P$ can be homotoped relatively to its boundary into $P$. We say the paring
is {\it acylindrical} if $\partial M\setminus P$ is incompressible and every incompressible
cylinder disjoint from $P$ and with boundary curves in $\partial M$ can be homotoped into $\partial M\setminus P$
relatively to $\partial M$.

We assume that we are given a word hyperbolic group $G$ as above and an embedding $\vp:\partial G\to S^2$, and we write  $\La_G= \vp(\partial G)$.
We will identify in the sequel subsets of $\partial G$ with subsets of $\La_G$ via the map $\vp$.

\REFPROP{prop:typeIIIplanar} Let $G$ be a one-ended hyperbolic group with planar boundary.
Let $v$ be a vertex of rigid type in its JSJ decomposition. 
Then the action of $G_v$ on $\La_v$ extends to a convergence
action of $S^2$ with limit set $\La_v$.
\ENDPROP

The proof of this proposition will require several steps, which we outline right now. We will first prove that 
we may find a degree $1$ map of $S^2$ transforming $\La_v$ onto a homeomorphic
copy of $Q_v$; this will enable us to prove that $Q_v$ is a degenerate carpet and that $G_v$ acts
as a geometrically finite convergence group (Prop.\,\ref{prop:rigide}). This action is planar and can be extended to a convergence action 
on $S^2$ with limit set the copy of $Q_v$ (Cor.\,\ref{lma:l4typeIII}).
We may then lift this action to an action of $S^2$ with limit set $\La_v$.

It will be used to prove:

\REFCOR{cor:typeIII} Let $G$ be a torsion-free one-ended hyperbolic group with a planar boundary and 
let $v$ be a vertex of rigid type in its JSJ decomposition
such that all the stabilizers of its incident edges are isomorphic to $\Z$. We also assume that $\confdim_{AR} G_v <2$.
Then there exists a compact hyperbolizable $3$-manifold $M_v$  with boundary and with fundamental group isomorphic to $G_v$
such that the conjugacy classes of incident edges define a maximal collection of incompressible disjoint simple closed curves
on $\partial M_v$.\ENDCOR
Thickening this multicurve into pairwise disjoint annuli provides us with an  acylindrical paring of $M_v$. 

We start by analysing the quotient action. The main step is provided by the next proposition.

\REFPROP{prop:rigide} Let $G$ be a one-ended hyperbolic group with planar boundary.
Let $v$ be a vertex of rigid type in its JSJ decomposition with stabilizer  $G_v$ and limit set  $\La_v$; the quotient $Q_v=\partial G/\sim_v$  is 
a degenerate carpet which satisfies the following properties.
\ben
\item There exists a pseudoisotopy $(\psi_t)_{t\in [0,1]}$ of the sphere such that,
writing $\psi=\psi_1$, $\psi(\La_v)$ 
is homeomorphic to $Q_v$ and such that fibers
are points except at non-trivial classes of $\sim_v$ where fibers are closed arcs. From now on, we let $Q_v=\psi(\La_v)$. 
\item The compact set  $\wLa_v=\psi^{-1}(Q_v)$ is connected and locally connected. 
\item  For any connected component $\oo$ of $S^2\setminus Q_v$, $Q_v\setminus \partial\oo$
is connected.
\item The action of $G_v$ on $Q_v$ is geometrically finite.
\een 
\ENDPROP

A pseudoisotopy is a continuous map $\psi:[0,1]\times S^2\to S^2$ where, for all $t\in[0,1)$, 
$\psi_t:x\mapsto \psi(t,x)$ is a homeomorphism of $S^2$
and $\psi_0$ is the identity. 

We start with a lemma which takes care of the points (1) and (2).

\REFLEM{lma:l1typeIII} There exists a pseudoisotopy $(\psi_t)_{t\in [0,1]}$ of the sphere such that,
writing $\psi=\psi_1$, $\psi(\La_v)$ 
is homeomorphic to $Q_v$ and such that fibers
are points except at non-trivial classes of $\sim_v$ where fibers are closed arcs. 
Moreover, $\wLa_v$ is connected and locally connected, where
$\wLa_v$ denotes the inverse image of $\psi(\La_v)$ under  $\psi$. \ENDLEM

\Proof 
Let $e\in E_v$ be an edge incident to $v$. Since $\La_e$ disconnects $Z_e$ from $\La_v$,
we may find an arc $c_e\subset Z_e$ joining $\La_e$.

Proceeding as above for all edges incident to $v$, we obtain a family of arcs $\{c_e\}_{e\in E_v}$. 
Since the sets $Z_e$ are disjoint,
these curves are pairwise disjoint. We wish to prove that the partition $\GGG$ of the sphere into these arcs and single points is upper semicontinuous, cf. \cite[Chap.\,VII]{whyburn:analytic_topology}. 
Since the non-trivial elements of this collection form a countable set, it is enough to prove that they form a null-sequence. This follows from the local connectivity of $\partial G$
since it implies that $(Z_e)_{e\in E_v}$ forms a null-sequence.

Since the elements of $\GGG$ are connected compact non-separating subsets of $S^2$, Moore's Theorem \cite[Theorem 25.1]{daverman:decompositions}
implies that the quotient $S^2/\GGG$ is homeomorphic to the sphere.  Note that $\GGG$ and $\sim_v$ agree on $\La_v$ so that
both quotients are homeomorphic (to $Q_v$). 
By \cite[Theorems 13.4, 25.1]{daverman:decompositions}, 
the decomposition $\GGG$ of $S^2$  has the property of being {\em strongly shrinkable}: 
there is  a pseudoisotopy $\psi_t: S^2 \to S^2, t \in [0,1]$ such that the fibers of $\psi$ agree with $\GGG$. 

The set $\wLa_v$ is clearly connected since $\psi$ is monotone. 
It remains to prove it is also locally connected. Let $x\in\wLa_v$ and let us consider a nested
family of connected neighborhoods $(V_n)$ of $\psi(x)$ with $\cap V_n=\{\psi(x)\}$. It follows that $\psi^{-1}(V_n)$
is also a sequence of nested connected neighborhoods of $x$ and that $\cap \psi^{-1}(V_n)=\psi^{-1}(\{\psi(x)\})$.
Therefore, we already know that $\wLa_v$ is locally connected at points which are fibers of $\psi$.
If $x$ belongs to the interior of some arc $c_e$, then we may also construct a basis of connected neighborhoods
as the arc is isolated in $\wLa_v$. If $x\in\La_e$, then $\psi^{-1}(V_n)\setminus c_e$ has 
two connected components, one of them ---$W_n$--- containing $x$ in its closure. We may then
add to $W_n$ a small subarc of $c_e$ to obtain a connected neighborhood $W_n'$ in $\wLa_v$ of $x$ 
so that $\cap W_n'=\{x\}$.
\endp


\Proof (Proposition \ref{prop:rigide}) 
 Note that $\La_v$ is one-dimensional (otherwise we would have $\La_v=S^2$), and $Q_v$ as well.
By Lemma \ref{lma:proptypeIII} and Lemma \ref{lma:degcarpet}, $Q_v$ is a degenerate carpet.
Points (1), (2) and (4) have been established in Lemmas \ref{lma:proptypeIII} and \ref{lma:l1typeIII}.

We now prove (3).
Let $\oo$ be a component of $S^2\setminus Q_v$, and let us prove that $Q_v\setminus\overline{\oo}$ is connected: let $Q_1$ and $Q_2$
be a partition of $Q_v\setminus\overline{\oo}$ into two open disjoint sets. Let $V$ be a connected component of $S^2\setminus (Q_v\cup\oo)$.

Note that $\partial{V}$ can intersect $\partial\oo$ at most at a single point, otherwise we would be able to disconnect $Q_v$ by removing two points,
contradicting Lemma \ref{lma:proptypeIII}. Therefore $\partial V\cap (Q_v\setminus\partial\oo)$ is connected: it is either contained in $Q_1$ or in $Q_2$. 

We may now define, for $j=1,2$, $U_j$ to be the union of $Q_j$ with all the components $V$ the boundaries of which are contained in $Q_j$: we obtain in this
way two disjoint open sets covering $S^2\setminus\overline{\oo}$. Since the latter is connected, $U_1$ or $U_2$ is empty, hence $Q_v\setminus\partial\oo$
is connected as well.
\endp

\REFCOR{lma:l4typeIII} The action of  $G_v$ on $Q_v$ can be extended to a convergence action on $S^2$ with limit set $Q_v$
such that $(S^2\setminus Q_v)/G_v$ is a finite union of surfaces of finite type.
\ENDCOR

\Proof According to Proposition \ref{prop:rigide}, the image of a boundary  component of $Q_v$ by any element of $G_v$  does not separate $Q_v$, so it bounds another boundary component,
hence the action of $G_v$ is planar (and geometrically finite). 

Let us prove that there are only finitely many orbits of such components. 

Let $\PPP=\{\psi^{-1}(\partial\oo)\cap\La_v,\ \oo\in\pi_0(S^2\setminus Q_v)\}$. Since the action on $Q_v$ is planar, 
the collection $\PPP$ is $G_v$-invariant. Note that each element is contained in the boundary of a connected component 
of $S^2\setminus \wLa_v$, one in each;
since $\wLa_v$ is a locally connected continuum
by Lemma \ref{lma:l1typeIII}, we conclude that $\PPP$ is  a $G_v$-invariant null-sequence.
Therefore, Proposition \ref{prop:ginvns} applies and proves that each stabilizer of $K\in\PPP$  is quasiconvex with limit set $K$
and that $\PPP$ is composed of finitely many orbits.
Pushing down the action of $G_v$  by $\psi$, we obtain finitely many orbits of boundary components, and 
the stabilizer $G_K$ of each $K\in\PPP$  
provides us with a geometrically finite action 
of the stabilizer $G_{\Omega}$ 
of  the boundary of each complementary component $\oo$  of $Q_v$.

 It follows from Theorem \ref{thm:extcva}  and Proposition \ref{prop:homeo}
that the action of $G_v$ to $Q_v$ is the restriction of a convergence action on $S^2$.

Finally, for each component $\Omega$ of $S^2\setminus Q_v$, the action of its stabilizer
$G_{\Omega}$ is geometrically finite, hence isomorphic to a geometrically finite Fuchsian group and
we may conclude that $\oo/G_v$ is a surface of finite type. 
\endp

\Proof (Proposition \ref{prop:typeIIIplanar})
We want to lift under $\psi$ the action of $G_v$ given by Corollary \ref{lma:l4typeIII}. 
Each parabolic point $p_e=\psi(\La_e)$ belongs to the boundary of two complementary components 
according to Lemma \ref{lma:proptypeIII}. 

Let  $H_e\subset (S^2\setminus Q_v)$ be the union of two closed $G_e$-invariant horocycles (one in each component)
that we choose small enough so that  their  $G_v$-orbit  are all pairwise disjoint.

Set $Z=\overline{\psi^{-1}(S^2\setminus (\cup_{e\in E_v} \psi^{-1}(H_e))}$.  
In other words, $Z$ is the whole sphere where we have taken off the union of the interiors of $\psi^{-1}(H_e)$, so it contains in particular $\La_v$. 
We first lift the action of $G_v$ to
$Z$. Let $g\in G_v$; if $x\in \La_v$, then we set $\hat{g}(x)=g(x)$; otherwise, we set $\hat{g}(x)=\psi^{-1}\circ g\circ \psi(x)$.
Let us prove that $\hat{g}$ is continuous. It is enough to consider a sequence $(x_n)$ disjoint from $\La_v$ which tends to some point $x$. 
For any limit point $z$ of $(\hat{g}(x_n))$, we have $\psi(z)= (g\circ \psi)(x)$.  The continuity of $\hat{g}$ follows at every $x$ in a trivial
fibre of $\psi$.  It remains to deal with points $x$ in the limit set of an edge group. 
Let $e\in E_v$ and let us consider a neighborhood $V_e$ of $p_e$ and set $W_e=(V_e\setminus H_e)\cup\{p_e\}$. If $V_e$ is suitably chosen, then $W_e\setminus\{p_e\}$
consists of two components $W_e^\pm$  which lift under $\psi$ to two disjoint open connected sets, each of which contains a single point of $\La_e\subset\La_v$ in its
closure. The same property holds for $g(W_e)$ and $\psi^{-1}(g(W_e))$. Let us assume that $x\in\La_e$ and that it is in the closure of $W_e^+$; then $\hat{g}(W_e^+)=\psi^{-1}\circ g\circ\psi(W_e^+)$
and   $\hat{g}(W_e^+)$ accumulates a unique point from $g(\La_{e})$.
Since $x$ is not isolated in $\La_v$, this point is $g(x)$. 
This shows that $\hat{g}$ is continuous, hence a homeomorphism. 


We now extend the action to the whole sphere.
Let $e$ represent an element of $E_v/G_v$; the set $S_e=\psi^{-1}(H_e)$ is  homeomorphic to a Jordan
domain and $G_e$ acts as an elementary convergence group on its boundary. 
Since  $E_v/G_v$ is finite,  Theorem \ref{thm:extcva}
implies that 
the action of $G_v$ extends to a convergence action on $S^2$.
\endp
 
We now prove the corollary.

\Proof (Corollary \ref{cor:typeIII}) According to Proposition \ref{prop:typeIIIplanar}, $G_v$ acts on $S^2$ with limit set $\La_v$.
Since its Ahlfors regular conformal dimension is strictly less than two, Corollary \ref{cor:planaraction} provides us with
a hyperbolizable manifold.

Each generator of a conjugacy class of an incident edge defines a non-trivial curve $\g_e$ in $M_v= \HH^3/G_v$.
Since the curves $\{c_e\}$ are pairwise disjoint off of $\La_v$, it follows that the curves $\g_e$ are homotopic to  
simple and pairwise disjoint curves on $\partial M_v$. If this family was not maximal, we could split $M_v$ 
along an incompressible annulus disjoint from the $\g_e$'s. But this would imply that the JSJ-decomposition of $G$ was not maximal. 
\endp

\subsection{Regular decomposition} 
We will focus on particular JSJ-decompositions which are suited to manifolds.

\begin{defn}[Regular JSJ decomposition] Let $G$ be a one-ended hyperbolic group
and let us consider its JSJ decomposition. We say it is regular if the following properties
hold:
\bit
\item[--] every two-ended group $H$ which appears as a vertex or an edge group is isomorphic
to $\Z$ and stabilizes the components of $\partial G\setminus \La_H$;
\item[--] the stabilizer of vertices of surface type is free;
\item[--] the stabilizer of vertices of rigid type is torsion-free.
\eit\end{defn}

The main point comes from the following proposition.

\REFPROP{prop:jsjreg} Let $G$ be a one-ended hyperbolic group with planar boundary and with 
a regular JSJ decomposition. We assume that the Ahlfors regular conformal dimension of each
vertex of rigid type is stricly less than $2$. Then $G$ is isomorphic to the fundamental group of a compact
hyperbolizable $3$-manifold with boundary.\ENDPROP

\Proof We first notice that since each elementary group $G_v$ or $G_e$ fixes the components of
$\partial G\setminus \La_v$ and $\partial G\setminus \La_e$, they are generated by primitive elements
of $G$.

If $v$ is a vertex of elementary type, we associate a solid torus $M_v= \SS^1\times \DD$  on which we consider on its
boundary pairwise disjoint incompressible annuli $\SS^1\times \al_e$ in bijection with its incident edges,
where $\al_e\subset\partial\DD$ are arcs.

If $v$ is of surface type, then $G_v$ is Fuchsian so it uniformizes a surface $S_v$;  we let  $M_v=S_v\times [0,1]$ and, 
noting that each incident edge in $T/G$ corresponds to a simple closed curve $\g_e\subset S_v$ bounding a hole
of $S_v$, we may associate on $\partial M_v$ disjoint annuli $\g_e\times [0,1]$. 

If $v$ is of rigid type, Corollary \ref{cor:typeIII} enables us to associate a pared manifold $M_v$.

Now, $T/G$ provides us with a manual to build a $3$-manifold $M$ with fundamental group isomorphic to $G$  
by gluing the $M_v$'s along the  annuli \cite{scott:wall}.
Thurston's hyperbolization theorem for Haken manifolds shows that $M$ is hyperbolizable. \endp



\section{The QCERF property}\label{sec:qcerf}
A group $G$ satisfies the  {\it  QCERF property} if every quasiconvex subgroup $A<G$ is {\it separable} i.e., 
for any $g\in G\setminus A$, 
there exists a finite index subroup of $G$
which contains $A$ but not $g$. This property will provide us with regular JSJ decompositions.


\REFTHM{thm:qcerf} Let $G$ be a non-elementary 
hyperbolic group which is QCERF.
There is a 
finite index torsion-free subgroup $H<G$ such that all the one-ended groups arising in its splitting
into a free product have regular JSJ decompositions.
\ENDTHM

The key point is the following group-theoretic result.

\REFPROP{prop:indexsbgroup} Let $A'<A<G$ be groups with $[A:A']<\infty$ and $A'$ separable in $G$.
Then there exist subgroups $A''$ and $H$ with the following properties:
\ben
\item $H$ is a normal subgroup of finite index in $G$;
\item $A''=H\cap A'$ is a normal subgroup of finite index in $A$ ;
\item for all $g\in G$, $(gAg^{-1})\cap H = gA''g^{-1}$.
\een \ENDPROP

\Proof Since $A'$ has finite index in $A$, we may pick a set 
of representatives $\{a_0,\ldots, a_n\}$ of $A/A'$ with $a_0$ the neutral element. Since $A'$ 
is separable, there is some finite index subgroup $G_j < G$ which
contains $A'$ but not $a_j$ for all $j\in\{1,\ldots, n\}$. 
Note that this implies that $a_j A'\cap G_j=\emty$.

Let $G'=\cap G_j$. Then $G'$ is a finite index subgroup of $G$ with
the property that $G'\cap A= A'$. Let $H$ be the largest  finite  index subgroup of $G'$
normal in $G$ and set $A''= A\cap H$.  

It follows that $H$ is a normal subgroup of finite index in $G$ and that
 $A''$ is a normal subgroup of finite index in $A$. Note that since $A\cap G'=A'$
we also have $A''=H\cap A'$.
Since $H$ is a normal subgroup, it follows that, for all $g\in G$,
$(gAg^{-1})\cap H = g(A\cap H)g^{-1}=gA''g^{-1}$.\endp

We include a proof of the following folklore result:
\REFCOR{cor:virtorsionfree} A hyperbolic group with the QCERF property is virtually torsion-free.\ENDCOR

\Proof A hyperbolic group has finitely many conjugacy classes of torsion elements \cite[Prop.\,4.13]{ghys:delaharpe:groupes}.
Let $\{g_1,\ldots, g_n\}$ be a set of representatives. Set $A_j$ to be the cyclic subgroup generated by $g_j$
and apply Proposition \ref{prop:indexsbgroup} to $\{e\}<A_j <G$ to obtain $G_j$ and set $G'=\cap G_j$.
We let the reader check that $G'$ is a finite index torsion-free subgroup of $G$.\endp

\Proof (Thm\,\ref{thm:qcerf}) Since $G$ is hyperbolic and QCERF, 
it contains a finite index subgroup $H'$ which is torsion-free (and QCERF) by Corollary \ref{cor:virtorsionfree}. 
It follows that its action on $\partial G$ is faithful. We first write $H'$ as a free product of one-ended groups
with a free group, cf. \S\,\ref{sec:strongacces}. Let $K'$ be one of the factors which is one-ended but not Fuchsian.

 Let
us consider a set of representatives of elementary vertex groups $\{K_v'\}$ arising from  its JSJ decomposition.
Each one of them contains a cyclic subgroup 
$A_v$ of finite index which stabilizes all the components of $\partial K\setminus \La_v$. 
Note that, for any $g\in G$, $gA_v g^{-1}$  satisfies the same conditions
at the vertex $g(v)$.

We now apply Proposition \ref{prop:indexsbgroup} to each triple $(A_v,K_v',H')$ 
to obtain finitely many finite index subgroups of $H'$ and we let $H$ denote their intersection.

By construction, $H$ is normal and of finite index in $H'$. 
Note that a virtually free torsion-free group is free \cite{stallings:free} so that we only need
to check the condition on the elementary type vertices and edges.
For a given factor   $K'$ as above,  we let $T_K$ denote the Bass-Serre tree given by its JSJ-decomposition.
Then $K=H\cap K'$ acts on $T_K$ and $T_K/K$ provides its JSJ-decomposition.
For each vertex $v$,  
we have $K_v=K_v'\cap H$.  Therefore, Proposition \ref{prop:indexsbgroup} ensures that $K_v$ 
is a subgroup of $A_v$. 
Moreover, each  edge group $K_e$ is cyclic since it is torsion-free and it fixes the components of $\partial K\setminus \La_e$ since 
two adjacent vertices are not of the same type and if it is incident to a vertex of  elementary type, then it already
fixes the complementary components.  This guarantees that $K$ admits a regular JSJ decomposition.
\endp

\section{Group actions on CAT(0) cube complexes}\label{sec:cccat}
A cube complex $X$ is a
CW-complex where each $n$-cell is a standard Euclidean $n$-cube and such that 
\ben
\item each closed cube is embedded into $X$;
\item the intersection of two cubes is either empty or a face.
\een
A cube complex is naturally endowed with a length structure such that each $n$-cell is isometric to
a unit Euclidean cube of the same dimension. We will focus on cube complexes which satisfy the CAT(0)
condition. We shall say that a group $G$ is {\it cubulated} if it admits a geometric and cellular action on a CAT(0) cube complex $X$.
We refer for instance to \cite{sageev:cc95} and \cite{bridson:haefliger:book} for details.

We gather some definitions and properties for future reference.

\subsection{Hyperplanes} \label{sec:hyperplane}
A fundamental feature of CAT(0) cube complexes comes from hyperplanes. Let us first define a {\it midcube} of a cube $[(-1/2),1/2]^n$
to be the intersection of the cube with a linear hyperplane orthogonal to one axis $[(-1/2),1/2]^n\cap\{x_j=0\}$, for some $j\in\{1,\ldots, n\}$.
A {\it hyperplane} is a maximal convex subset of a CAT(0) cube complex for which the intersection with any cube is either empty
or a midcube.  

Hyperplanes have many interesting properties \cite{sageev:cc95,haglund:finiteindex}. Among them:
\ben
\item given an edge, there is a unique hyperplane which intersects it orthogonally; 
\item  a hyperplane divides a CAT(0) cube complex into exactly two connected components which are both convex.
\een

Say a hyperplane $Y\subset X$ is {\it essential} if, for any $R>0$, none of the two components of $X\setminus Y$ is contained in 
the $R$-neighborhood of $Y$. Note that if $X$ is proper and hyperbolic, then 
any pair of distinct points $x,y\in\partial X$ is separated by an essential hyperplane.
If we also assume that $X$ supports a geometric action of a group $G$, then
the action of the stabilizer $H$ of a hyperplane $Y$ is always cocompact. 

\subsection{Special actions}\label{sec:special}
Haglund and Wise have defined a particular class of non-positively curved cube complexes named {\it special}, see \cite{haglund:wise:special}  for the precise
definition and a proper introduction to the subject. They enjoy two properties which will be of interest for the present work.

Say a hyperbolic group $G$ has a {\it special action} if it acts  cellularly and geometrically 
on a  CAT(0) cube complex $X$ such that $X/G$ is special.
In this case,
\ben
\item  the group $G$ splits over the stabilizer of any essential hyperplane; 
\item the group $G$ has the QCERF property  \cite[Thm\,1.3]{haglund:wise:special}. 
\een
A group is {\it virtually special} if it contains a finite index subgroup which admits a special action on a CAT(0) cube complex.
By \cite[Thm\,1.3 and Lma 7.5]{haglund:wise:special}, a virtually special group has also the QCERF property.

Following Wise \cite{wise:qcvxh,agol:virtualHaken}, 
define the class $\QQQ\VVV\HHH$ as the smallest class of hyperbolic
groups that contains the trivial group $\{1\}$ and is closed under the following operations:
\bit
\item[--] if $G=A\star_C B$ with $A,B\in\QQQ\VVV\HHH$ and $C$ f.g. and quasiconvex in $G$ then $G\in\QQQ\VVV\HHH$;
\item[--] if $G=A\star_C$ with $A\in\QQQ\VVV\HHH$ and $C$ f.g. and quasiconvex in $G$ then $G\in\QQQ\VVV\HHH$;
\item[--] if $H < G$ with $H\in\QQQ\VVV\HHH$ and $[G:H]<\infty$ then $G\in\QQQ\VVV\HHH$.
\eit
A group in $\QQQ\VVV\HHH$ is said to have a {\it quasiconvex virtual hierarchy}.

 Wise  proved that  a hyperbolic group has a quasiconvex virtual hierarchy if and only if it is virtually special 
\ \cite[Thm.\,13.1]{wise:qcvxh}. Moreover, Agol proved that any cubulated hyperbolic group
admits a finite index subgroup with a special action \cite[Thm.\,1.1]{agol:virtualHaken}.
In summary, we have

\REFTHM{thm:agolwise} Let $G$ be a hyperbolic group. The following are equivalent
\bit
\item[--] $G$ is cubulated;
\item[--] $G$ is virtually special;
\item[--] $G$ admits a quasiconvex virtual hierarchy.
\eit
In this situation, $G$ satisfies the QCERF property and there is a finite index subgroup of $G$ which admits a quasiconvex hierarchy. \ENDTHM

We record the following application.

\REFCOR{prop:virtspeck} If $G$ is a convex-cocompact Kleinian group, then $G$ is virtually special.\ENDCOR

\Proof By Brooks' theorem \cite{brooks:extension}, we may assume that it is a quasiconvex subroup of a cocompact
Kleinian group. By \cite{bergeron:wise:cubulation} and  Theorem \ref{thm:agolwise}, the latter is virtually special so
\cite[Thm\,H]{haglund:finiteindex} implies that $G$ is as well.\endp

\subsection{Strong accessibility}\label{sec:strongacces}
What follows is inspired by similar properties of compact $3$-manifolds. If $M$ is a compact $3$-manifold with  compressible boundary, it can be cut along {\it compression disks} into finitely many pieces
each of which has incompressible boundary. The fundamental group is split as a free product, with factors given
by the fundamental groups of each piece. Given a compact hyperbolizable manifold $M$ with incompressible
boundary, we may cut it into finitely many pieces along essential annuli so that the remaining pieces
are acylindrical pared manifolds; an annulus  $A$ is {\it essential} if  it is properly embedded, incompressible and {\it non-peripheral} i.e,  the inclusion $i:A\to M$
is not homotopic to a map $f:A\to M$ such that $f(A)\subset\partial M$. This corresponds to the JSJ decomposition of the group given in \S\,\ref{sec:jsj}.

Let $G$ be a non-elementary hyperbolic group. If it is not one-ended then it splits over a finite group \cite{stallings:yale}.
By \cite{dunwoody:accessibility}, there is a quasiconvex splitting over finite groups such that each vertex group is finite
or
one-ended; when $G$ is torsion-free, it leads to a free product of a free group with finitely many
one-ended groups.  We may then consider the JSJ-decomposition of the remaining one-ended vertex subgroups and proceed inductively.
If $G$ has no element of order two, then this process stops in finite time 
\cite{delzant:potyagailo:access, louder:touikan}, exactly as for Haken manifolds which are cut into finitely many balls. 

In the end, in the group setting,
we are left with finite groups, virtually Fuchsian groups and/or one-ended hyperbolic groups with no local cut points 
in their boundaries. 
If $G$ has planar boundary, then those latter groups are carpet groups \cite[Thm\,4]{kapovich:kleiner:lowdim}. 

\subsection{A criterion for the QCERF property}
This section is devoted to the proof of the following proposition.

\REFPROP{thm:virtuallyspe} Let $G$ be a non-elementary hyperbolic group with a planar
boundary different from the sphere. The group $G$ is QCERF if one of the following conditions hold:
(a)  $\confdim_{AR} G<2$; (b)   $\confdim_{AR} H<2$ for every subgroup arising as a rigid vertex group of the JSJ decomposition of $G$ or of one of its maximal one-ended subgroups; 
(c) $G$ has no elements of order two and if $\confdim_{AR}H<2$ 
for all carpet quasiconvex subgroups $H$,
then $G$ is QCERF. \ENDPROP

\Proof 
According to Theorem \ref{thm:agolwise}, it suffices to prove that $G$ has quasiconvex virtual hierarchy to conclude
that $G$ has the QCERF property.

 Note that the assumption   (a) implies  (b) by \cite[Prop.\,2.2.11]{mackay:tyson:confdim}, so let us assume  that (b) holds. 
Let us split $G$ over finite groups so that
every factor is elementary or one-ended. Then $G$ will have a quasiconvex virtual hierarchy provided
every rigid group arising in the JSJ-decomposition is also in $\QQQ\VVV\HHH$.
But such a group has a planar action by Proposition \ref{prop:typeIIIplanar} so Theorem \ref{thm:main} implies
that it is virtually Kleinian under the assumptions (a) or (b). So Corollary \ref{prop:virtspeck} enables us to conclude that $G$ is virtually special.

We now assume that $G$ has no elements of order $2$. 
From the strong accessibility of such groups (\S\,\ref{sec:strongacces}), 
we just have to deal with carpet groups. Assuming their conformal dimension is
strictly less than two, Corollary \ref{cor:main} implies that they are virtually isomorphic to convex-cocompact Kleinian
groups so are virtually special according to Corollary  \ref{prop:virtspeck}.
\endp


\section{Characterizations of convex-cocompact Kleinian groups}
Let $G$ be a word hyperbolic group with a planar boundary. 

If $G$ is elementary, virtually free or if $\partial G$ is homeomorphic to the unit circle,
then it is already known that $G$ contains a finite-index subgroup isomorphic to a convex-cocompact
Kleinian group.
\ben
\item If $G$ is finite, then we may consider the trivial group.
\item If $G$ is two-ended, then it contains a finite index subgroup isomorphic to $\lag z \mapsto 2 z\rag$.
\item If $G$ is virtually free, then  there is a finite index free subgroup $H$ which is isomorphic to the fundamental
group of a handlebody. This implies that $H$ is isomorphic to a convex-cocompact Kleinian group ---a so-called Schottky group.
\item If $\partial G$ is homeomorphic to the unit circle, then Theorem \ref{thm:qmonS1} implies that  $G$ contains a finite index subgroup
isomorphic to a cocompact Fuchsian group. 
\een

Thus, we assume throughout this section that $G$ is not in one of the above classes.
As in \S\,\ref{sec:strongacces}, we may first
split $G$ over finite groups so that each vertex group is finite
or
one-ended. Let $\HHH_1$ denote each of those vertex groups, (set $\HHH_1=\{G\}$ if $G$ is one-ended)
and let us denote by $\HHH_1'$ those vertex groups in $\HHH_1$ which
are one-ended and non-virtually Fuchsian. For each one of them, we may consider its JSJ decomposition as in \S\,\ref{sec:jsj}.

In this section, the term {\it rigid group} will refer to the stabilizer of a rigid vertex coming from the JSJ-decomposition of a group
in $\HHH_1'$.

\subsection{Dynamical characterization of groups with planar boundaries}

We first recall that a free group is the fundamental group of a handlebody or a solid torus if its rank is one,
so it can always be uniformized by a convex-cocompact Kleinian group.

We start with a proposition which sums up the previous sections:

\REFPROP{prop:basechar} Let $G$ be a non-elementary 
hyperbolic group with a planar boundary. If the Ahlfors-regular conformal dimension of all rigid groups 
is strictly less than two, then the group $G$  is virtually isomorphic to a convex-cocompact Kleinian group.\ENDPROP

\Proof
According to Proposition \ref{thm:virtuallyspe},
$G$ is QCERF and so Theorem \ref{thm:qcerf} implies that it contains a torsion-free subgroup $H$ of
finite index such that the JSJ-decompositions of its one-ended free factors are regular.

We construct a compact Haken manifold with fundamental group isomorphic to $H$ as follows.

Let us first write $H=H_0\star H_1\star \ldots \star H_n$ as a free product of a free group $H_0$ and of one-ended subgroups.
Now, each non-elementary $H_j$ has conformal dimension strictly less than $2$, since there are finite-index subgroups of groups with
the same property (implying that their gauges coincide). 
Either $H_j$ is free or with boundary isomorphic to the circle or $H_j$ satisfies the assumptions
of Proposition \ref{prop:jsjreg}. In either case, $H_j$ is isomorphic to 
the fundamental group of a hyperbolizable manifold $M_j$ with boundary. We may then glue these manifolds
along disks on their boundaries to obtain a Haken  manifold $M$ with fundamental group isomorphic to $H$.
Theorem \ref{thm:thurston} proves that $M$ is hyperbolizable, so that, up to index two, $H$ is isomorphic
to a convex-cocompact Kleinian group.\endp

We may now establish  Theorem \ref{thm:main1}. 

\Proof (Theorem \ref{thm:main1}) The necessity part  is due to Sullivan, see
\cite{sullivan:ihes50}. 

Let $G$ be a non-elementary hyperbolic group with planar
boundary of Ahlfors regular conformal dimension strictly less than two. According to \cite[Prop.\,2.2.11]{mackay:tyson:confdim},
the Ahlfors-regular conformal dimension of every rigid group has also Ahlfors regular conformal dimension strictly less than two.
Therefore, Proposition \ref{prop:basechar} applies.\endp




Let us deduce Theorem \ref{cor:main1}: 

\Proof (Theorem \ref{cor:main1}) Let $\vp:\partial G\to \cbar$ be a quasisymmetric embedding, and set $\La=\vp(\partial G)$.
It follows from Proposition \ref{prop:basechar} 
that it is enough to prove that the rigid  groups 
have Ahlfors regular conformal
dimension strictly less than $2$. We may as well assume $\La$  to be connected.

So let $v$ be such a vertex with stabilizer $G_v$,  limit set $\La_v\subset\La$
and boundary components $\PPP_v$. 
Note that the set  $\La_v$ together with its complementary components in $\cbar$
define a separation structure in the sense of \S\,\ref{sec:geoprop}. Moreover, the action of $G_v$ on $\La_v$  is planar by Proposition \ref{prop:typeIIIplanar}. 
Therefore
\ben
\item the set $\La_v$ is relatively porous to $\PPP_v$ according to Proposition \ref{prop:ginvns};
\item  the boundary components are uniform quasicircles by Proposition \ref {prop:sierpgroup}, so that the complementary
components are uniformly porous.
\een
Hence, Proposition \ref{prop:porosity} implies that $\La_v$ is porous in $\cbar$. This is equivalent to having 
its Assouad dimension strictly less than two \cite[Thm.\,5.2]{luukkainen:assouad}. This in turn
implies that $G_v$ has Ahlfors-regular conformal dimension strictly less than two \cite[Prop.\,2.2.11]{mackay:tyson:confdim}. \endp


We now turn to  Theorem \ref{thm:main0}; we just prove (3) implies (1).

\Proof (Theorem \ref{thm:main0}) 
Let $G$  be a non-elementary hyperbolic group with planar boundary, no elements of order $2$ and with carpet quasiconvex
subgroups of conformal dimension strictly less than two. Proposition  \ref{thm:virtuallyspe} ensures that
$G$ is QCERF. Therefore, we may  assume that $G$ is torsion free by Corollary \ref{cor:virtorsionfree}. 
We will argue by induction on the length of the hierarchy provided
by its strong accessibility. We let $\HHH_1$ denote the set of vertices arising when splitting into freely indecomposable subgroups, and $\HHH_1'\subset\HHH_1$
 be the  non-Fuchsian one-ended factors. Define 
$\HHH_2$ to be the  non-Fuchsian rigid vertex groups  arising
in the JSJ splittings of the elements in $\HHH_1'$. We proceed
inductively so that $\HHH_{2n+1}'$ denotes the factors of the decomposition as a free
product by a free group with one-ended groups of the elements of $\HHH_{2n}$ and $\HHH_{2n+2}$ consists
of the vertices 
obtained by the JSJ-decomposition of the non-Fuchsian one-ended elements of $\HHH_{2n+1}$. 
The strong accessibility ensures that this procedure stops in finitely many steps \cite{vavritchek:access, louder:touikan}.

 
We will prove that $H$ is isomorphic to the fundamental group of a hyperbolizable manifold by induction on its hierarchy.
The initial cases $n=1,2$ follow from Proposition \ref{prop:basechar}. 

Assuming Theorem \ref{thm:main0} holds up to rank $2n$, we show how to deal with a group $H$ with $2n+2$
generations. Note that each element of $\HHH_2$ has rank at most $2n$, so they are virtually isomorphic to convex-compact Kleinian groups. 
It follows that their Ahlfors regular conformal dimension
is strictly less than $2$. Therefore, Proposition \ref{prop:basechar} applies again and completes the proof.
\endp

\subsection{Cubulated groups with planar boundaries}
We will deduce Theorem \ref{thm:main2} from the following:

\REFTHM{thm:cubplanar} Let $G$ be a quasiconvex subgroup of a hyperbolic cubulated group
$\hat{G}$ with boundary homeomorphic to the two-sphere. Then $G$ is virtually a convex-cocompact Kleinian group.
\ENDTHM

We start by preparing the ambient group $\hat{G}$:

\REFPROP{prop:cubspherespe} Let $G$ be a cubulated hyperbolic group with 
 boundary homeomorphic to the two-sphere. There exists a finite index torsion-free subgoup which
 admits a special action on a cube complex such that every hyperplane  stabilizer is isomorphic
 to a cocompact Fuchsian group. \ENDPROP
 
 \Proof  Let $X$ be a CAT(0) cube complex on which $G$ acts geometrically. 
 
Let $x,y\in\partial X$ be two distinct points. We first prove that we may separate them with the limit set of 
a cocompact Fuchsian subgroup of $G$. We start with a hyperplane of $X$ with stabilizer $H$ which separates $x$ and $y$. 
Considering a quasiconvex subgroup, we may assume that
 $H$ is one-ended. It follows that $x$ and $y$ lie  in two different components of $\partial X\setminus \La_H$.
Note that the action of $H$ on $\partial X$ is planar.  Therefore, by Proposition \ref{prop:sierpgroup}, there is a Fuchsian
group which separates $x$ and $y$.

By \cite{bergeron:wise:cubulation}, $G$ acts on  a CAT(0) cube complex $Y$ such that the stabilizer of each hyperplane
is virtually Fuchsian. By Theorem \ref{thm:agolwise} and Corollary \ref{cor:virtorsionfree}, 
there is a finite index torsion-free subgroup of $G$ which has a special action
on $X$. The stabilizers of the hyperplanes in  this subgroup are now Fuchsian.\endp

 \Proof (Thm \ref{thm:cubplanar}) The proof  will proceed by induction: 
 we may start with $G<\hat{G}$ where
 $\hat{G}$ is torsion-free and admits a special action on a CAT(0) cube complex $\hat{X}$, cf.
 Proposition \ref{prop:cubspherespe}.  Since $G$ is quasiconvex, its action on $X$ is convex-cocompact according
 to \cite[Thm.\,H]{haglund:finiteindex}: there is a convex subcomplex $X\subset \hat{X}$ invariant by $G$ with
 $X/G$ is compact. 
 
 Let $Y\subset X$ 
 be an essential hyperplane of $X$. Set $\YYY=\cup_{g\in G} g(Y)$. We define a graph
 $T$ as follows: the vertices are the connected components of $X\setminus \YYY$ and two vertices form
 an edge if they are separated by exactly one hyperplane $g(Y)$ for some $g\in G$. Since 
 $G$ has a special action, it follows that $T$ is a tree and that the action of $G$ on $T$ is simplicial, minimal and
 without edge inversions, cf. \cite{haglund:wise:special}. If we let $C$ be the stabilizer of $Y$, then we have shown
 that $G$ is either an amalgamated product $G=A \star_C B$ or an HNN extension $G=A \star_C$, where $A$ and $B$ are
 stabilizers of components $X\setminus\YYY$: their action are also convex-cocompact. 

If no vertex group contains a carpet
 group, then Corollary \ref{cor:main0} shows that these vertex groups are Kleinian, hence their
 Ahlfors regular conformal dimension is strictly less than two. We then apply the following proposition, the
 proof of which is postponed later on:
 
 \REFPROP{prop:induction} With the notation above, if the Ahlfors regular conformal dimension of every vertex
 group is strictly less than two, then $G$ is conjugate to a convex-cocompact Kleinian group. \ENDPROP
 
 This ends the proof in this case.
 Otherwise, each vertex group is convex-cocompact in $X$ so we may find a convex subcomplex invariant by the vertex groups, and
 proceed as above until their vertex groups contain no carpet group.  
 This process ends since the action of $G$ is special i.e., $G$ admits a quasiconvex hierarchy.
We obtain a rooted finite tree of quasiconvex subgroups where children of a vertex correspond to a splitting of
 it.  So, we may apply, as above, Corollary \ref{cor:main0}  to the leaves and inductively Proposition \ref{prop:induction}
 in order to reconstruct the whole group $G$. 
 \endp
 
 \Proof (Prop.\,\ref{prop:induction}) We first treat the case $G=A\star_C B$. 
 It follows from Corollary \ref{cor:planaraction} that $A$ and $B$ are conjugate to convex-cocompact Kleinian groups.
 
 The hyperplane $Y$ is defined by an orthogonal edge $e\in X^{(1)}$, cf. \S\,\ref{sec:hyperplane}. This is also
 an edge of $\hat{X}$ so it defines a hyperplane $\hat{Y}\subset\hat{X}$. It follows that
 $Y=\hat{Y}\cap X$ so $\stab_G Y= \stab_{\hat{G}}\hat{Y}\cap G$ and $C$ is a quasiconvex subgroup of a cocompact
 Fuchsian group $\hat{C}$. Moreover, $\hat{Y}$ is clearly inessential with respect to $A$ and $B$ so  we may name
 the connected components of 
 $\partial\hat{X}\setminus \La_{\hat{C}}$ $D_A$ and $D_B$ so that $D_A\cap \La_A=\emty$ and
 $D_B\cap \La_B=\emty$. Moreover, it follows from \cite[Cor.\,4.6]{martin:tukia:msri} that  the action of $C$ on $\partial\hat{X}$
 is globally conjugate to a convex-cocompact Fuchsian group. Hence
 there is an equivariant involution $\iota_C:\partial \hat{X}\to\partial \hat{X}$ which fixes $\La_{\hat{C}}$ pointwise and exchanges
 its complementary components $D_A$ and $D_B$.
 
 It follows that $(\overline{D_A}\setminus \La_C)/C$ is a subsurface of a boundary component of the Kleinian manifold $M_A$,
 $(\overline{D_B}\setminus \La_C)/C$ is also  a subsurface of a boundary component of the Kleinian manifold $M_B$ and 
 $\iota_C$ induces an orientation reversing homeomorphism between them.
 We may then define
 $M_G= M_A \sqcup_{\iota_C} M_B$: this is  a Haken manifold with fundamental group isomorphic to $G$
 \cite{scott:wall}. Therefore, Theorem \ref{thm:thurston} implies that $G$ is Kleinian and Corollary \ref{cor:planaraction} that
 the actions are conjugate.
 
 If $G= A\star_C$, the proof is similar: there is some element $g_0\in G\setminus A$ such that the HNN extension
 is obtained by identifying $C$ with $C'=g_0 C g_0^{-1}$. Since the action is special, $C$ and $C'$ will be contained
 in cocompact Fuchsian groups $\hat{C}$ and $\hat{C}'$ which bound disjoint disks $D_C$ and  $D_{C'}$ disjoint from $\La_A$.
 We may then 
 glue the compact surface $(\overline{D_C}\setminus\La_C)/C$ with $(\overline{D_{C'}}\setminus\La_{C'})/C'$, both contained in the boundary of $M_A$, to obtain a Haken
 manifold with fundamental group isomorphic to $G$.
  \endp

 We now prove Theorem \ref{thm:main2}.
 
 \Proof (Thm\,\ref{thm:main2}) The necessity comes from Corollary \ref{prop:virtspeck}.
 For the suffciency, Theorem \ref{thm:main0} tells us that we just need to deal with groups
 with boundary a carpet or a sphere. For the latter, Theorem \ref{thm:cubplanar} shows that
the group is virtually Kleinian.  

For the carpet case, let $G$ be a carpet group and let $H_1$, \ldots, $H_k$ denote representatives of the peripheral Fuchsian groups,
cf. Proposition \ref{prop:sierpgroup}. Let us take another copy $(G',H_1',\ldots, ,H_k')$ of $(G,H_1,\ldots, ,H_k)$ and consider
the graph of groups with vertices $G$ and $G'$ and with $k$ edges identifying each $H_j$ with $H_j'$.  
According to  \cite[Thm\,5]{kapovich:kleiner:lowdim}, one obtains a hyperbolic group $\hat{G}$ with boundary
the sphere so that $G$ has become a quasiconvex subgroup of $\hat{G}$. 

Since $G$ is cubulated, it follows from Theorem \ref{thm:agolwise} that $G$ admits a quasiconvex hierarchy, so $\hat{G}$ as well
and we may conclude that $\hat{G}$ is cubulated.


We may now apply Theorem \ref{thm:cubplanar} and conclude that $G$ is a virtually convex-compact Kleinian group.
 \endp
 
 

 \subsection{Quasi-isometric rigidity}\label{sec:qirigidity}
 
 We prove Theorem  \ref{thm:cmain}, since Theorem  \ref{thm:mmain} is a particular case of it.  
 Let $G$ admit  a quasi-isometric embedding into $\HH^3$. Then its image $Z$ is quasiconvex, hence 
 $G$ is word hyperbolic, and we may endow its boundary with a metric from its gauge. 
The quasi-isometry extends into  a quasiymmetric homeomorphism between $\partial G$ and $\La_Z$  (see Remark \ref{rmk:qivsqs}) 
so Theorem \ref{cor:main1} implies that $G$ is virtually a convex cocompact Kleinian group. \endp


We take this opportunity to provide some additional consequences for carpet groups which seem relevant to the present discussion, and
which are direct corollaries of \cite{bonk:kleiner:merenkov}, but which seem to
have been overlooked in the literature; see nonetheless  \cite{kapovich:kleiner:lowdim, frigerio:mostow,frigerio:commensurability,bonk:icm:qcgeom,kleiner:icm2006}
 for similar properties. Important and recent progress on such questions, especially on the rigidity of carpet groups, can be found in \cite{merenkov:carpetgrouprigidity:preprint}.

\begin{proposition} Let $G$ be a carpet group quasi-isometric to a Kleinian group
and which acts faithfully on its boundary. 
Then $G$ has finite index in the group $G_M$ of quasi-M\"obius homeomorphisms of $\partial G$. 
The following properties also hold.
\ben
\item The group $G_M$ is the unique maximal word hyperbolic group 
in the quasi-isometry class of $G$ which acts faithfully on its boundary. It is isomorphic to a Kleinian group.
\item Whenever $G$ acts geometrically on a proper geodesic metric space, any self-quasi-isometry of $X$ lies at bounded distance from
an element of $G_M$.
\item For any  group $H$ quasi-isometric to $G$, $H/F$  
is isomorphic to a finite index subgroup of $G_M$ where $F$ is the kernel of the action of $H$ on $\partial H$. 
\item For any convex-cocompact Kleinian group $K$ quasi-isometric to $G$, there is a finite
locally isometric covering $p:M_K\to\HH^3/G_M$. 
\item Any group quasi-isometric to $G$ is commensurable to $G$. 
\een
\end{proposition}
For point (2), since $X$ is quasi-isometric to $G_M$, the latter admits a quasi-action by quasi-isometries on $X$.
Point (5) was already known, cf. \cite{kapovich:kleiner:lowdim,frigerio:commensurability}. It is also known that in general
Kleinian groups which are quasi-isometric need not be commensurable \cite{whyte:amenability}.
Let us note that  a torsion-free convex-cocompact Kleinian group is a carpet group exactly when it uniformizes an {\it acylindrical compact manifold} i.e., a manifold
with an  incompressible boundary and no essential annuli.


\Proof Since $G$ is quasi-isometric to a Kleinian group, it is a virtually convex-cocompact Kleinian group. 
According to Remark \ref{rmk:gdcboundary}, there is a quasisymmetric homeomorphism of $\partial G$ to a round carpet $\La$ (the complementary components are spherical disks) with spherical measure zero. Therefore, any quasi-M\"obius self-homeomorphism of $\La$ is the restriction of a M\"obius transformation 
according to \cite{bonk:kleiner:merenkov}.

This implies that the group of quasi-M\"obius homeomorphisms $G_M$  of $\La$ is a subgroup of the group of M\"obius transformations of the sphere. 
It is clearly discrete since any sequence which tends uniformly to the identity will have to eventually fix at least  three circles, 
implying that such a sequence is eventually the identity. Therefore, $G_M$ is a convergence group, and since it contains $G$, its action is uniform on $\La$, 
hence it is word hyperbolic. Since the boundaries of $G$ and $G_M$ coincide, $G$ has finite
index in $G_M$ \cite{kapovich:short:greenberg}. 
The group $G_M$ is clearly maximal since any group quasi-isometric to $G$ acts on $\La$ by M\"obius transformations
as well with a finite kernel. 

The trace at infinity of any quasi-isometry is a quasi-M\"obius homeomorphism, hence its trace belongs to $G_M$, and this
implies it is at bounded distance from the corrresponding element.

Any (M\"obius)  group $K$ quasi-isometric to $G$ which acts faithfully on $\La$ is a subgroup of $G_M$.
Therefore, any such group has finite index in $G_M$ so is commensurable to $G$.
\endp

To end this paper, besides the Cannon and Kapovich-Kleiner conjectures, the quasi-isometric rigidity of convex cocompact groups motivates the following questions, 
see also \cite[Question 12.6]{kapovich:kleinmore}:

\noindent{\bf Questions.---} 
\ben
\item  {\em Is the class of all finitely generated Kleinian subgroups of $\PP SL(2,\C)$ quasi-isometrically rigid?}
\item {\em  Is the class of fundamental groups of compact $3$-manifolds quasi-isometrically rigid?}
\een

Thanks to Thurston's hyperbolization and the work above, 
the  first problem reduces to the case of groups which are geometrically finite, without rank 1 maximal parabolic subgroups, cf. \cite[Thm. 19.6]{kapovich:book}. 
In the case of lattices in $\PP SL(2,\C)$, the positive answer is  due to Sullivan and Schwartz. 
For the second question,  it is known that the class of irreducible manifolds of zero Euler characteristic  is rigid. This follows from
works of Perel'man, Sullivan, Schwartz, Kapovich and Leeb, Gromov, Pansu, and Eskin, Fisher and Whyte.


\bibliographystyle{math}
\bibliography{refs}

\end{document}